\documentclass{article}

\usepackage{tabularx} 
\usepackage{amsmath}  
\usepackage{graphicx} 
\usepackage[margin=1in,letterpaper]{geometry} 
\usepackage{cite} 
\usepackage[colorlinks,allcolors=black,pdftex]{hyperref}
\usepackage{amsfonts}
\def\vec#1{\ensuremath{\mathbf{#1}}}
\usepackage[position=top]{subfig}
\usepackage{pgfplots}
\usepackage[title]{appendix}
\usepackage{booktabs}
\usepackage[capitalise]{cleveref}
\usepackage{amssymb}
\usepackage{authblk}

\usepackage{placeins}
\providecommand{\keywords}[1]{\textbf{\textit{Keywords}} #1}
\pgfplotsset{compat=newest}

\hypersetup{
	colorlinks=true,       
	linkcolor=blue,        
	citecolor=blue,        
	filecolor=magenta,     
	urlcolor=blue
}
\usepackage{blindtext}

\begin{document}

\title{Modelling of physical systems with a Hopf bifurcation using mechanistic models and machine learning}
\author[1]{{K.~H.~ Lee}\thanks{Corresponding Author: jz18526@bristol.ac.uk}}
\author[1]{{D.~A.~W.~Barton}\thanks{david.barton@bristol.ac.uk}}
\author[2]{{L.~Renson}\thanks{l.renson@imperial.ac.uk}}
\affil[1]{Department of Engineering Mathematics, University of Bristol, UK}
\affil[2]{Department of Mechanical Engineering, Imperial College London, UK}
\date{}
\maketitle

\begin{abstract}
We propose a new hybrid modelling approach that combines a mechanistic model with a machine-learnt model to predict the limit cycle oscillations of physical systems with a Hopf bifurcation. The mechanistic model is an ordinary differential equation normal-form model capturing the bifurcation structure of the system. A data-driven mapping from this model to the experimental observations is then identified based on experimental data using machine learning techniques. The proposed method is first demonstrated numerically on a Van der Pol oscillator and a three-degree-of-freedom aeroelastic model. It is then applied to model the behaviour of a physical aeroelastic structure exhibiting limit cycle oscillations during wind tunnel tests. The method is shown to be general, data-efficient and to offer good accuracy without any prior knowledge about the system other than its bifurcation structure.

\end{abstract}


\section{Introduction}
Limit-cycle oscillations (LCOs) are periodic responses that can be observed in many systems such as aircraft wings~\cite{dimitriadis2017introduction}, wheels~\cite{beregi2019bifurcation}, machine tools~\cite{kalmar2001subcritical}, and living cells~\cite{adimy2005stability,guo2003hopf}.
\textit{Self-excited} systems are a common source of LCOs and are typically modelled using ordinary differential equations (ODEs) where the variation of a parameter beyond a critical value (a bifurcation point) triggers oscillations. 

Deriving a low-dimensional mathematical model that quantitatively captures the onset and amplitude of LCOs is usually a very challenging problem as self-excited systems are typically characterised by the interplay of several physical phenomena. Take, for example, the aforementioned fluid-structure and tyre-ground interactions in wings and wheels, respectively.

For self-excited systems with a Hopf bifurcation, the change in parameter leads to a loss of stability of the equilibrium and the birth of a family of LCOs near the bifurcation point. To analyse the periodic responses of such systems and determine the amplitude of the LCOs, it is customary to reduce the dynamics of the original system near the bifurcation point to a low-dimensional centre manifold. The reduced system is topologically equivalent to the full model, i.e. the vector flows are locally qualitatively identical and there exists a smooth invertible coordinate transformation between them~\cite{kuznetsov2013elements}. The hybrid modelling methodology proposed here takes advantage of this equivalence. A low-dimensional normal form model is used to capture the phenomenology of the real system, i.e. its bifurcation structure. A data-driven mapping from this model to the full system is then identified using machine learning (ML) and experimentally measured data (here LCOs). 

The use of ML models is attractive because they are theoretically able to represent any continuous functions~\cite{winkler2017performance}. However, using traditional ML techniques and models presents a number of difficulties. For instance, a considerable amount of data is often needed to train ML models~\cite{wang2021understanding}. Furthermore, even if such data is available and the obtained ML model accurately represents the training data set, ML models can still fail to generalise to unseen conditions~\cite{rasmussen2003gaussian} or even fail to capture the fundamental physics of the system~\cite{kim2019dpm}. The approach proposed in this paper is inspired by the recent development of scientific machine learning (SciML), which aims at making ML models more interpretable, more consistent with the known laws of physics, and less data-hungry by combining them with mechanistic (i.e. physics-based) models. Examples of SciML approaches include Physics-Informed Neural Networks (PINNs) \cite{raissi2019physics} where a neural network is used to solve and discover partial differential equations (PDEs) while respecting the laws of physics through constraints incorporated into the training cost function. In Refs.~\cite{raissi2017machine,raissi2017inferring} linear ordinary differential equations (ODEs) were successfully discovered from data using probabilistic machine learning and Gaussian process regression. Universal differential equations (UDEs)~\cite{rackauckas2020universal} are differential equation models combining mechanistic differential equations with universal approximators such as neural networks, Chebyshev expansions, or random forests directly introduced into the model equations. 

While existing studies have mostly focused on numerical simulations and a quantitative agreement between time-series at particular parameter values, the approach proposed here aims to capture the bifurcation diagram of a physical system, which requires the accurate prediction of the system's parameter dependence, its long-term behaviour (here LCOs), and to deal with experimental (i.e. noise-corrupted) measurements. In this context, Beregi et al. \cite{beregi2021using} combined machine-learnable functions with mechanistic models to capture bifurcation diagrams. The approach proposed here differs in that it uses only the knowledge of the bifurcation structure observed in the experiment. The normal form model and its associated bifurcation structure form the mechanistic model that captures the ``\textit{physics}'' of the system and underpins the otherwise data-driven model. The hybrid mechanistic/machine-learnt (M/ML) model obtained does not rely on any problem specific variables and is therefore applicable to any system exhibiting the bifurcation structure imposed by the underlying mechanistic model, even without any other physical model available. 

Simple parameter sweeps are often enough to reveal such a bifurcation structure, including its supercritical or subcritical nature. The use of a model to capture the type of bifurcation observed in the data is very beneficial as it reduces the amount of data required to train the model and improves the ability of the model to interpolate between data points and even extrapolate outside the range of control parameters used for model training. Moreover, by leveraging the fact that the dynamics of the system evolves on a low-dimensional sub-manifold and using a polar representation, the training of time series requires the time integration of a single one-dimensional ODE, which is computationally much more efficient than other data-driven modelling procedures using numerical integration of the entire model as in \cite{rackauckas2020universal}. The hybrid M/ML models developed in this paper could be exploited in different ways. For instance, the machine-learnt mapping could be analysed to improve understanding about the physical system and provide new insights into the derivation of more accurate mechanistic models. The proposed models have also the potential to be used as digital twins where the nature of the system is captured qualitatively using the mechanistic part of the model and data is continuously used to refine and evolve the model during the system's life. Exploring these applications is outside the scope of the paper. 

The manuscript is structured as follows. The second section discusses the theoretical background of the proposed modelling approach, and the third section explains the process used to train the data-driven part of the model. In Section 4, the proposed modelling approach is validated numerically on a Van der Pol oscillator and a three-degree-of-freedom aeroelastic model representative of the physical system studied in Section 5. The numerical experiments show that the proposed methodology accurately predicts the phase portrait, bifurcation diagram and time series of the studied systems. Finally, the proposed hybrid modelling approach is successfully validated experimentally in Section 5. The physical system considered is an aerofoil that exhibited LCOs when tested in the wind tunnel. The data exploited in this paper was collected using control-based continuation (CBC) \cite{sieber2008control,renson2019application,renson2016robust,barton2017control,barton2013systematic}. Contrary to conventional parameter sweeps which can only measure stable LCOs, CBC uses feedback control to stabilise, and hence observe unstable LCOs, which can benefit parameter estimation~\cite{deCesare22,beregi2020improving}. The hybrid M/ML model built with the proposed approach is shown to capture the overall bifurcation structure of the physical system and to reproduce quantitatively the amplitude of both stable and unstable LCOs measured in the experiment. 

\section{Model structure}\label{sec:ModelStructure}

\paragraph{Mechanistic model} It is assumed that the physical system of interest can be represented by an unknown continuous-time dynamical system~$(N,\Phi^t)$, where $N \subset \mathbb{R}^{n+1}$ with $n \geq 2$ is the number of states and $\Phi^t: \mathbb{R}^+ \times N \rightarrow N$ is the evolution of the flow governed by a set of ODEs. The system $(N,\Phi^t)$ is assumed to undergo a Hopf bifurcation at certain parameter $\mu=\mu_0$ and the sub-/super-critical nature of the bifurcation is known (or at least identifiable from experimental data). The parameter $\mu$ is constant over time and taken as one of the dimensions of $N$. 

There exists a 3-dimensional centre manifold $M^c$ near the Hopf bifurcation point that can be parametrised as a graph $\vec M^c_{\vec{x}}$ \cite{carr2012applications}:
\begin{align}\label{center_manifold}
\vec M^c_{\vec{x}}=\{(\vec x,\vec y)| \vec y=\vec h(\vec x)\},
 \end{align}
\noindent where $\vec x$ is the center subspace, $\vec y$ is the attracting subspace defined at the equilibrium of $(N,\Phi_{\text{full}}^t)$, and $\vec h$ is a nonlinear function. $M^c$ is an attracting invariant manifold in the state-space of the full dynamical system~$(N,\Phi_{\text{full}}^t)$. The dynamics of the system on the centre manifold is defined as $(M^c,\Phi_{\text{red}}^t)$, where $\Phi_{\text{red}}^t: \mathbb{R}^+ \times M^c \rightarrow M^c$ is the evolution of the flow on the centre manifold $M^c$. The system $(M^c,\Phi_{\text{red}}^t)$ is a reduced order model of $(N,\Phi_{\text{full}}^t)$ with $\text{dim}(M^c) \leq \text{dim}(N)$. The reduced dynamics is usually expressed as $\dot {\vec x}=\vec f (\vec x)$ where $\vec f$ is the system's vector field projected onto the centre manifold parametrised by $\vec x$. It is also possible to find a smooth, invertible change of coordinates such that the reduced dynamics can be represented using the modified Hopf normal form
\begin{align}\label{NF-model}
  \begin{split}
  \dot{u}_1 &= (\mu-\mu_0) u_1-u_2\Omega(u_1,u_2,\mu)+a_2 u_1(u_1^2+u_2^2)-u_1(u_1^2+u_2^2)^2,\\
  \dot{u}_2 &= (\mu-\mu_0) u_2+u_1\Omega(u_1,u_2,\mu)+a_2 u_2(u_1^2+u_2^2)-u_2(u_1^2+u_2^2)^2,\\
  \dot{\mu} &=0,
  \end{split}
 \end{align}

\noindent which can also be written in polar coordinates as
\begin{align}\label{NF-model-polar}
  \begin{split}
  \dot{r} &= (\mu-\mu_0) r+a_2 r^3 - r^5\\
  \dot{\theta} &= \Omega(r,\theta,\mu)\\
  \dot{\mu} &=0.
  \end{split}
\end{align}
\noindent where $(u_1, u_2, \mu)$, or $(r, \theta, \mu)$, are the coordinates parameterizing the invariant manifold. The sign of the coefficient $a_2$ depends on the criticality of the Hopf bifurcation. Fifth-order terms are added to the standard Hopf normal form to introduce a saddle-node bifurcation of periodic orbits and capture the presence of stable LCOs frequently observed in systems with subcritical Hopf bifurcations. For the latter, stable and unstable LCOs coexist for $\mu \in [\mu_0-a_2^2/4, \; \mu_0]$. $\Omega$ is an a priori unknown function that governs the speed of the oscillations and that
reproduces the time evolution of the data (see \cref{sec:ModelTraining}). From \cref{NF-model-polar}, it is clear that the oscillation amplitude $r$ is independent of the oscillation speed. This observation will be exploited in Section~\ref{sec:ModelTraining} to simplify the ML model training by learning the coordinate mapping and the oscillation speed separately.

\cref{NF-model,NF-model-polar} are not the only way to represent the dynamics of a system with a Hopf bifurcation followed by a saddle node bifurcation. An alternative parameterization of the invariant manifold, using, for instance, $\vec x$ as in \cref{center_manifold}, can be used. This would result in a set of ODEs that is different from \cref{NF-model,NF-model-polar}. However, those different reduced systems are topologically equivalent to each other, i.e. one can be transformed to another by a change of coordinates~\cite{kuznetsov2013elements}. The choice to take the well-known normal form as a mechanistic model was made to emphasize the nature of the phenomenon targeted by the model.

\begin{figure}[tb]
	\centering
  \includegraphics[width=0.8\linewidth]{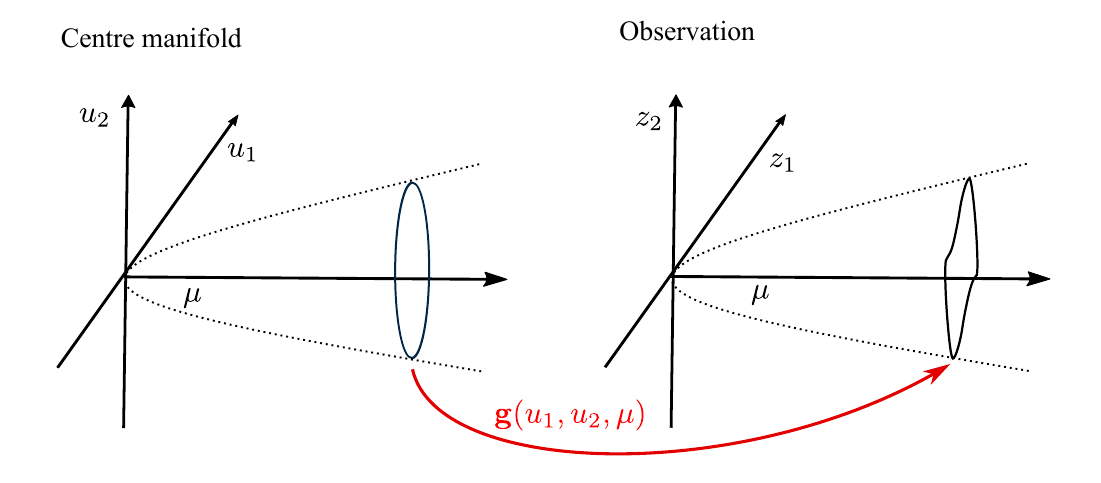}
  \caption{Geometric illustration of the coordinate transformation $\vec U_{12}(u_1, u_2, \mu)$ applied to a LCO of a supercritical Hopf bifurcation.}
  \label{Utrans}
\end{figure}

\paragraph{Mapping to observations} Following the definition of the mechanistic model~\eqref{NF-model}, a data-driven transformation from the model to the observation is defined. Let's consider the measured observations $\vec z=[z_1,\ldots,z_m,\mu]^\text{T} \in \mathbb{R}^{m+1}$, where $m$ is the total number of states observed, and the predicted observations $\hat {\vec z}=[\hat z_1,\ldots,\hat z_m,\hat \mu]^\text{T} \in \mathbb{R}^{m+1}$. A function $\vec g(u_1,u_2,\mu) = {\vec{\hat{z}}}$ can be defined to map the dynamics of the reduced system~\eqref{NF-model} to the predicted experimental observations $\hat {\vec z}$ made on the centre manifold $M^c$. The map $\vec g$ can be defined as a vector of two functions $\vec g= \left[ \vec U^T, \; g_\mu \right]^T $. The first function, $\vec U(u_1,u_2,\mu) = [\hat z_1,\ldots,\hat z_m]^\text{T}$, represents the mapping between $(u_1, u_2,\mu)$ and the observed states. The LCOs in Eq.~\eqref{NF-model} trace circular trajectories in the plane $(u_1, u_2)$. The objective of the mapping $\vec U$ is thus to transform these circles into the distorted closed curves observed experimentally (as illustrated in \cref{Utrans}). The second part of the map $\vec g$, $g_\mu(u_1,u_2,\mu) = \hat \mu$, represents the mapping between the model parameter $\mu$ and the predicted observed parameter. In the present and common case where the bifurcation parameter is directly measured during the experiment and not re-scaled, the mapping reduces to a simple projection, i.e. $g_\mu(u_1,u_2,\mu)=\mu$.

\section{Model training}\label{sec:ModelTraining}

\subsection{Closed orbit representation}\label{sec:shape}

We assume there exist a pair of coordinates for which the measured LCOs form closed curves that can be parameterised in polar coordinates, i.e. the LCOs form curves that do not self-intersect and have a unique angular parameterisation. This assumption is satisfied in the neighbourhood of the Hopf bifurcation point, and it is assumed that it extends to all the measured LCOs. For convenience, this pair of coordinates is labelled $(z_1,z_2)$, while the remaining measured coordinates are $(z_3, ..., z_m)$. The map $\vec U$ is split accordingly as $\vec U= \left[ \vec U_{12}^T, \; \vec U_{3\ldots m}^T \right]^T$. 

The map $\vec U_{12}$ from the normal form coordinates $(u_1, u_2, \mu)$ to the predicted observations $(\hat{z}_1, \hat{z}_2)$ is first sought. The particular challenge associated with finding this first map $\vec U_{12}$ is that the correspondence between points $(u_1,u_2)$ in the normal form coordinates and observations $(z_1,z_2)$ is initially unknown. Therefore, it is not possible to obtain $\vec U_{12}$ by solving a regression problem as input and output data points cannot be paired together. However, once this first map is found, the rest of the map $\vec U$ addressing the presence of additional states $(\hat{z}_3, \dots, \hat{z}_m)$ can be determined easily. This will be discussed later in this section. 

\paragraph{Finding the first mapping} To train $\vec U_{12}$, the idea is to compare continuous representation of the predicted and measured LCOs. This approach has the advantage of avoiding any point-wise comparison between data points and model predictions. The training process starts by taking a user-defined number of points along periodic responses in the normal form coordinates. Those points are mapped to the observation space using the current estimate of the mapping $\vec U_{12}$. Following the coordinate transformation, the closed curves obtained from the transformed trajectories can be directly compared with the measured LCOs. Comparing closed planar curves is a well-established problem in pattern recognition~\cite{rosenfeld1976digital,zhang2004review}, and a popular way to approach this problem is to use a Fourier representation of the curve along the arc-length~\cite{zahn1972fourier}. However, with such a representation, it is difficult to define a metric between two distinct curves if they do not share at least one point. For this reason, we here consider the simpler approach of directly using a phase-like angle to parameterise the orbit. This assumption is consistent with the normal form model and the closed curves observed experimentally, that have a much simpler geometry than the one usually investigated in pattern recognition~\cite{zahn1972fourier}. The polar representation of the planar orbits in terms of amplitude and angle is obtained for measured and predicted curves as
\begin{equation}\label{R_theta}
  R=\sqrt{z_1^2+z_2^2}, \quad \theta=\tan^{-1} \frac{z_2}{z_1}, \quad \quad \text{and} \quad \quad
  \hat{R}=\sqrt{\hat{z}_1^2+\hat{z}_2^2}, \quad \hat{\theta}=\tan^{-1} \frac{\hat{z}_2}{\hat{z}_1}.
\end{equation}
\noindent As the polar representations of the LCOs are assumed to be smooth and periodic, they can be represented as a truncated Fourier series
\begin{align}\label{FS}
  R(\theta)=a_0 + \sum_{k=1}^{n_h} a_k \cos({k\theta}) + \sum_{k=1}^{n_h} b_k \sin({k\theta}), \quad \quad \text{and} \quad \quad \hat{R}(\hat{\theta})=\hat{a}_0 + \sum_{k=1}^{n_h} \hat{a}_k \cos({k\hat{\theta}}) + \sum_{k=1}^{n_h} \hat{b}_k \sin({k\hat{\theta}}),
\end{align}
\noindent where the number of Fourier modes, $n_h$, is assumed to be sufficiently large to have a small approximation error. The shape of the closed curves $R(\theta)$ and $\hat{R}(\hat{\theta})$, are thus represented by vectors of coefficients determined by the Fourier projection $\Phi(\cdot)$ as
$\Phi(R):C_p([0,2\pi], \mathbb{R})	\to \mathbb{R}^{2n_h+1} =
[a_0,\ldots,a_{n_h},b_1,\ldots,b_{n_h}]^\text{T}$ where
\begin{equation}\label{Galerkin}
    a_0 = \frac{1}{2\pi}\int_0^{2\pi} R\mathrm{d}\theta, \qquad a_n = \frac{1}{ \pi}\int_0^{2\pi} R\cos(n\theta)\mathrm{d}\theta, \qquad b_n = \frac{1}{ \pi}\int_0^{2\pi} R\sin(n\theta)\mathrm{d}\theta \qquad \text{for} \: n=1,2,\dots,n_h.
\end{equation}
\noindent Alternatively, the vector of coefficient can be computed in a least-square sense directly using \cref{FS} and the pseudo inverse \cite{penrose1955generalized}.

Taking a family of LCOs from the branch emerging at the Hopf bifurcation point, the error between model predictions and the data is given by 
\begin{align}\label{XiU}
  \Xi_{\vec U}=\sum_{i=1}^{m_s} \| \Phi(R_i)-\Phi(\hat R_i) \|
\end{align}
\noindent where $m_s$ is the number of measured LCOs. This measure of the model error is the cost function that is minimised during the training of the coordinate transformation detailed in Section~\ref{sec:map_func}. 

\paragraph{Mapping to other coordinates} After obtaining $\vec U_{12}$, it is possible to use the inverse transformation $\vec U_{12}^{-1}$ to find the points in the normal form coordinates that are associated with the measured data. Considering those points as inputs, it is then straightforward to train the $m-2$ remaining maps to the output observations $(z_3, \dots, z_m)$ using standard input-output regression techniques such as kernel ridge regression \cite{kanagawa2018gaussian} or neural networks.

Although not explored here, an alternative approach to predict the remaining measured coordinates would be to consider the pair $(z_1, z_2)$ as so-called `\textit{master}' coordinates and $(z_3, \dots, z_m)$ as `\textit{slave}' coordinates of the system. Once the master coordinates are obtained from the normal form model using $\vec U_{12}$, a second mapping from the master coordinates to the slave coordinates can be learnt. This second problem is also a regression problem which can be easily solved. This procedure is conceptually similar to the classical centre manifold reduction approach (see \cref{center_manifold}). However, here, the mappings to the master coordinates and between master-slave coordinates are sought based on experimental data. This approach is also completely independent of the bifurcation structure. 

More generally, there exist also other approaches to train the overall mapping $\vec{U}$. For instance, a common parameterisation for the predicted (model) and observed (data) time series could be enforced. Time cannot be used as the oscillation speed $\Omega$ is trained separately. It is therefore natural to resort to a geometrical parameterisation of the time series in terms of phase-like angles like the one used for the parametrisation of the closed curves. Following a re-parametrisation of the predicted (measured) responses in terms of this phase angle, a one-to-one correspondence between data points in the normal form and observation coordinates can be assumed, and input-output regression performed to identify the mapping associated with each measured coordinate. While this approach appears straightforward, it was found more difficult to generate model predictions at the particular phase angles observed experimentally. As such, interpolation of the model predictions was necessary, which affected the precision of the overall training procedure.

\subsection{Functional form of the map}\label{sec:map_func}
Neural networks provide a flexible approach to model the mapping $\vec U_{12}$~\cite{lin2018resnet,winkler2017performance}. However, it was found that using a neural network alone often leads to mappings that do not preserve the topology of the LCOs. It is therefore advantageous to use a simpler initial transformation that preserves this topology. This transformation has the additional benefits of reducing the complexity of the neural network and simplifying its training. The planar mapping $\vec U_{12}(u_1,u_2,\mu) = [\hat z_1,\hat z_2]^\text{T} $ is thus defined as the sum of three separate contributions as

\begin{align}\label{g}
  \vec U_{12} (u_1,u_2,\mu)= \underbrace{T_L (u_1,u_2,\mu) \quad + \quad T_s}_{\text{Transforms LCO to an ellipse}}+\underbrace{\textrm{NN}_{\Theta_{\vec U}} (u_1,u_2,\mu)}_{\text{Correction}},
  \end{align}
\noindent where
\begin{equation}
    T_L (u_1,u_2,\mu) = \begin{bmatrix}
      l_{11} & l_{12} & l_{13} \\
      l_{21} & l_{22} & l_{23} \\
      l_{31} & l_{32} & l_{33}
     
  \end{bmatrix} 
  \begin{bmatrix}
    u_1 \\
    u_2 \\
    \mu
\end{bmatrix} \quad \quad \text{and} \quad \quad 
  T_s =   
  \begin{bmatrix}
    s_1 \\
    s_2 
\end{bmatrix}.
\end{equation}
The linear transformation performed by $T_L(u_1,u_2,\mu) = [\hat z_1, \hat z_2]^\text{T} $ stretches and rotates the closed curves. This transformation includes nine unknown parameters $l_{ij}$ that will be estimated using experimental data. The transformation matrix must be nonsingular. The coordinate transformation $T_s(u_1,u_2,\mu)= [\hat z_1, \hat z_2]^\text{T} $  applies a rigid translation of the trajectories and requires two additional parameters $s_1,s_2$. Finally, $\textrm{NN}_{\Theta_{\vec U}}(u_1,u_2,\mu) = [\hat z_1, \hat z_2]^\text{T} $ is a neural network with unknown weight vector $\Theta_{\vec U}$.
 
\begin{figure}[t!]
 \centering
 \includegraphics[width=1.0\linewidth]{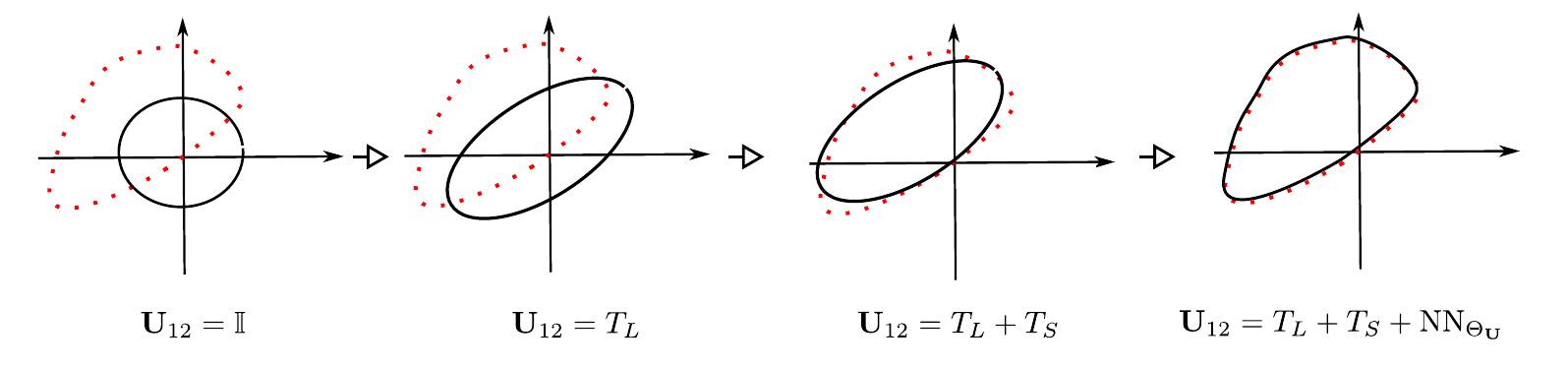}
\caption{Combined effect of the different components included in the transformation from normal form coordinates $\vec u$ to physical space $\vec z$. \textcolor{black}{($\boldsymbol{-}$)} LCO obtained after transformation. (\textcolor{red}{$\boldsymbol{\bullet}$}) Experimental data. }
 \label{T1}
\end{figure}

Including $T_L$ and $T_s$ explicitly in $\vec U_{12}$ can be interpreted as introducing additional physics or knowledge into the definition of the coordinate transformation. Indeed, shifting and rescaling coordinates must be performed to transform the orbits from the space of the normal form model to the space of the physical system and this operation is commonly performed in normal form calculations in bifurcation analysis. With this approach, the neural network model, $\textrm{NN}_{\Theta_\vec U}$, can be viewed as a ``small'' correction to the initial transformation performed by $T_L + T_s$ (\cref{T1}). 

\subsection{Oscillation speed}

The flow velocity on a LCO trajectory is defined by the general expression
\begin{align}\label{eq:Om_V}
\Omega(r\cos{\theta},r\sin{\theta},\mu)=\omega_0+\textrm{NN}_{\Theta_{\Omega}}(r,\mu)^\textrm{T} [1, \cos{\theta},\ldots,\cos{n_h \theta},\sin{\theta},\ldots,\sin{n_h \theta}]^\text{T},
\end{align}
\noindent where the first term, $\omega_0$, represents the fundamental oscillation frequency. The second term provides a periodic correction to this fundamental frequency in order to capture state and parameter dependencies. While this general correction term is particularly useful for capturing systems where multiple timescales occur within the LCO, it can be significantly simplified for systems where only one frequency dominates the response (see \cref{sec:NumVal}).

Following the training of the coordinate transformation, $\Omega$ is trained by minimising  the prediction error between predicted and observed time series. To generate time series from the model, suitable initial conditions in the normal form coordinates that corresponds to the initial measured data must be obtained. One approach to find the initial conditions $\vec u^{i}(t_1)$ would be to solve
\begin{align} \label{init1}
[\vec u^{i}(t_1),\mu_i]-\vec g^{-1}(\vec z^i(t_1),\mu_i)=0
\end{align}
\noindent where $t_1$ represents the first time instant in the time series. However, model inaccuracies and measurement noise perturbs the initial point, $\vec u^{i}(t_1)$, away from the trajectory of the LCO predicted at the parameter value $\mu_i$. While this is not a significant problem when training stable LCOs; it becomes an issue for unstable LCOs as the numerical integration of the initial value problem will not approach the trajectory of LCOs. The approach followed to solve this issue is to find the intial conditions $\vec u^i(t_1)$ for which the model prediction $(\hat z^i_1,\hat z^i_2) = \vec U_{12}(\vec u^i(t_1), \mu_i) $ has the same phase angle $\theta$ as the initial conditions of the measured signal. The initial conditions in the normal form coordinates, $\vec u^{i}(t_1)$, are thus found by solving 
\begin{align}\label{zp_ang}
\text {ang} \circ \vec U_{12}(\vec u^{i}(t_1),\mu_i)-\text {ang} \circ \vec z^i(t_1)=0,
\end{align}
\noindent where $\text{ang}: (x,y) \mapsto \tan^{-1}(y/x)$ measures the phase angle of the vector $[x,y]^T$.  \cref{zp_ang} can be solved using Newton method. Note that the initial conditions in the normal form coordinates cannot be set directly to the phase angle found in the data because the mapping $\vec U_{12}$ does not necessarily preserve this angle.

Once the initial conditions have been determined, the polar form of the normal form model is considered for the numerical integration. In this case, only the second equation of \cref{NF-model-polar} needs to be integrated as the LCOs correspond to fixed points of the first equation. Only integrating the second equation has also the advantage of avoiding any numerical instability issues even on the unstable solutions. Indeed, this second equation corresponds to the direction of the velocity vector, i.e. the direction vector of the trivial Floquet multiplier that is equal to unity.

The parameters $\omega_0$ and $\Theta_{\Omega}$ defining the oscillation speed $\Omega$ are then determined by minimizing the cost function 
\begin{align}\label{XiO}
\Xi_{\Omega}=\sum_{i=1}^{m_s} \sum_{j} \| \vec U_{12}(\vec u^{i}(t_j),\mu_i)-\vec z_i(t_j)  \|,
\end{align}

\subsection{Learning stages}\label{sec:learning_stages}
The model training is a three-stage process. The parameters of the linear transformations $T_L$ and $T_s$ are found first by minimising $\Xi_U$. An approximate value of the bifurcation parameter value $\mu_0$ is used during that process. After the training of the linear transformation, the LCOs in the normal-form space can be mapped to ellipses that are ``close'' to the measured trajectories. During the second training step, the parameters of $\textrm{NN}_{\Theta_{\vec U}}$ and more precise values for $\mu_0$ and $a_2$ are found by further minimising $\Xi_U$. The linear transformation parameters are kept constant during this process. The third training step is to find the parameters associated with the oscillation speed $\Omega$, i.e. $\omega_0$ and $\Theta_{\Omega}$, by minimizing $\Xi_{\Omega}$.

Traditional deep learning packages such as PyTorch\cite{paszke2019pytorch} and Flux.jl\cite{innes2018flux} can be used to train $\vec U$ and the other model parameters using optimisation techniques such as stochastic gradient decent method~\cite{da2014method,liu1989limited}. For $\Omega$, the package DiffEqFlux.jl \cite{rackauckas2020universal} which uses stochastic gradient descent methods on the solutions of differential equations~\cite{zhang2017discrete,lauss2018discrete} was used. 

\section{Numerical demonstration}\label{sec:NumVal}
In this section, the method developed in \cref{sec:ModelStructure,sec:ModelTraining} is demonstrated numerically on a Van der Pol oscillator and a 3-degree-of-freedom model of an aerofoil undergoing aeroelastic oscillations. The synthetic data used for model training is noise-free and was obtained using time integration. Demonstration on real experimental data is carried out in \cref{sec:ExpVal}.

\subsection{Van der Pol oscillator}\label{sec:vdp}
The equations governing the dynamics of the Van der Pol oscillator are
\begin{align}\label{vdp}
	\begin{split}
\frac{dz_1}{dt}&=z_2,\\
\frac{dz_2}{dt}&=2\mu z_2-z_1^2 z_2-z_1,
\end{split}
\end{align}
\noindent where the states $(z_1,z_2)$ and the control parameter $\mu$ are all assumed to be measured directly. For this example, a supercritical Hopf bifurcation occurs at $\mu_0 = 0$ and only stable LCOs exist. As such, the parameter $a_2$ of the mechanistic model~\eqref{NF-model} is set equal to $-1$ and the fifth-order terms are removed. Training data is generated for six different parameter values $\mu=(0.1,\;0.28,\; 0.46,\; 0.64,\; 0.82,\; 0.1)$. At each parameter values, the oscillator response is simulated over 10 seconds using initial conditions on the LCOs (i.e. there are no transient in the data) and a sampling time of 0.02 $s$. This represents 500 samples per time series, and hence 3000 samples for the whole training data set.

\begin{figure}[hbt!]
    \centering
    \subfloat{\includegraphics[width=0.7\textwidth]{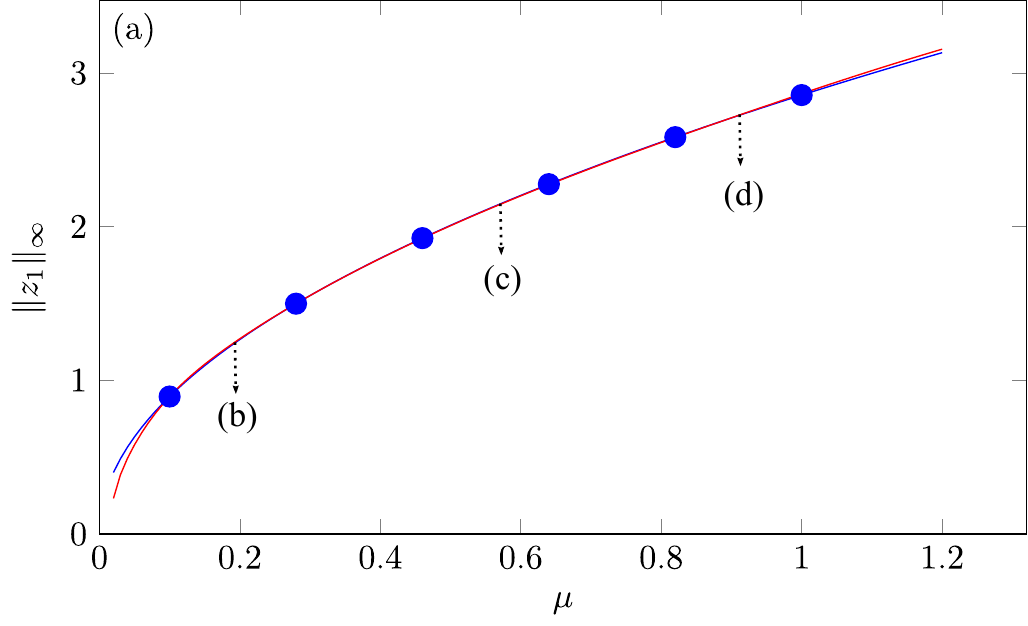}}
    \vspace*{-0.4em}
    \subfloat{\includegraphics[width=1.0\textwidth]{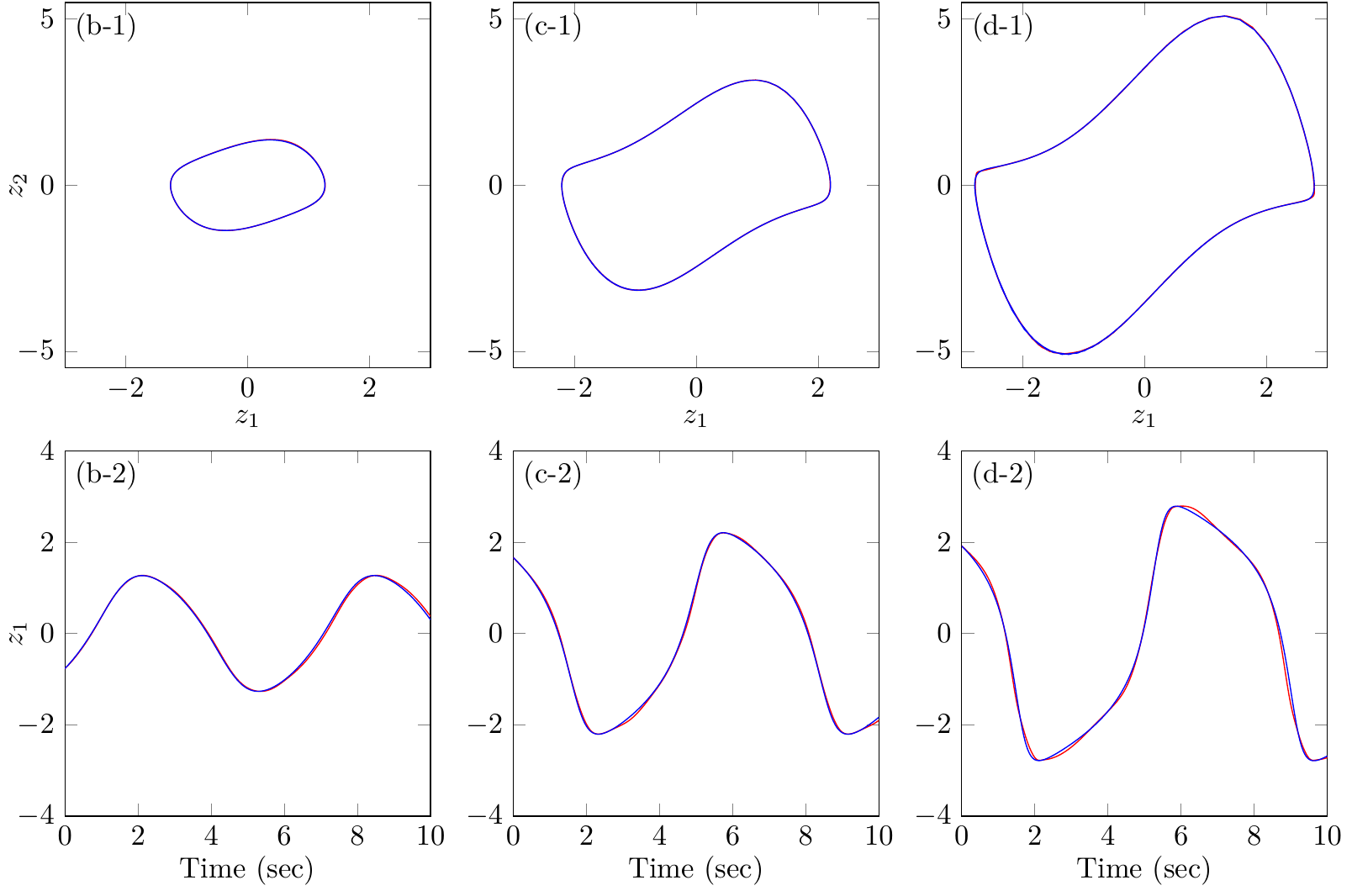}}
    \caption{Comparison between the Van der Pol model (\textcolor{blue}{$\boldsymbol{-}$}) and the hybrid M/ML model (\textcolor{red}{$\boldsymbol{-}$}). (a) Bifurcation diagram where (\textcolor{blue}{$\boldsymbol{\bullet}$}) are the LCOs used for model training. (b-1--d-1) Phase portraits and (b-2--d-2) time-series at the untrained locations reported on the bifurcation diagram.}
    \label{fig:vdp_all}
  \end{figure}

Following the procedure outlined in \cref{sec:ModelTraining}, the coordinate transformation is trained first by minimising $\Xi_{\vec U}$. The NN used within $\vec U$ consists of three inputs, two hidden layers each with 32 neurons and hyperbolic tangent activation functions, and two outputs. 300 iterations in ADAM~\cite{ruder2016overview} with a learning rate of $0.01$ followed by 1000 BFGS \cite{fletcher2013practical} iterations with a $10^{-5}$ learning rate were necessary to estimate the NN parameters $\Theta_{\vec U}$. A comparison between the bifurcation diagrams of the real and identified models shows that the hybrid M/ML model accurately captures the system's topological features (\cref{fig:vdp_all}(a)). Phase portraits are shown in Figs.~\ref{fig:vdp_all}(b-1)--(d-1). They further demonstrate that an accurate transformation from the normal-form coordinates to the physical coordinates is achieved for the range of parameter values considered. The parameter value at the bifurcation point was estimated at $\mu=0.02$, which is very close to the actual value $0$.

Following the training of the coordinate transformation, $\Omega$ is estimated. The $\textrm{NN}_{\Theta_{\Omega}}$ in model~\cref{eq:Om_V} is set to include three inputs, two hidden layers each with 32 neurons and hyperbolic tangent activation functions, and 13 outputs ($n_h=10$). The model parameters are estimated by minimizing $\Xi_{\Omega}$. The training was performed using 2000 iterations in ADAM with a learning rate of $0.01$ followed by 1000 BFGS iterations with a $10^{-5}$ learning rate were necessary. Figs.~\ref{fig:vdp_all}(b-2)--(d-2) show a very good agreement between the time series of the reference and identified models. As the bifurcation parameter increases, the time scale separation becomes more pronounced and errors become noticeable in the transition between the fast and slow portions of the time series (\cref{fig:vdp_all}(d-2)). Similar observations can be made for the other state (not shown for conciseness). 

\subsection{Aeroelastic model}\label{sec:NumVal_aero}
A 3-DOF aeroelastic system~\cite{abdelkefi2013analytical} is now considered to demonstrate numerically the proposed method. This model is qualitatively representative of the physical system tested in \cref{sec:ExpVal}. The equations of motion of this system are 
\begin{equation}\label{eq:flutter2}
  \vec{M} \ddot{\vec{x}} + \vec{D} \dot{\vec{x}} + \vec{K} \vec{x} + \vec{N}(\alpha) = 0,
\end{equation}
\noindent where 
\begin{subequations}
\begin{gather}\label{eq:2-3}
    \vec M =
      \begin{bmatrix}
          m_T+\pi \rho b^2            & m_w x_\alpha b-a\pi\rho b^3   & 0 \\
          m_w x_\alpha b-a\pi\rho b^3 & I_\alpha+\pi(1/8+a^2)\rho b^4 & 0 \\
          0                           & 0                             & 1
      \end{bmatrix},\\
    \vec D =
    \begin{bmatrix}
          c_h+2\pi\rho bU\hat{c}       & (1+\hat{c}(1-2a))\pi\rho b^2U                 & 2\pi U^2b(c_1c_2+c_3c_4) \\
          -2\pi(a+1/2)\rho b^2\hat{c}U & c_\alpha+(1/2-a)(1-\hat{c}(1+2a))\pi\rho b^3U & -2\pi\rho b^2U^2(a+1/2)(c_1c_2+c_3c_4) \\
          -1/b                               & a-1/2                                               & (c_2+c_4)U/b
    \end{bmatrix},\\
    \vec K =
    \begin{bmatrix}
          k_h & 2\pi\rho bU^2\hat{c}                  & 2\pi U^3c_2c_4(c_1+c_3) \\
          0   & k_\alpha-2\pi(1/2+a)\rho\hat{c}b^2U^2 & -2\pi\rho bU^3(a+1/2)c_2c_4(c_1+c_3) \\
          0   & -U/b                                        & c_2c_4U^2/b^2
    \end{bmatrix},
  \end{gather}
\end{subequations}
\noindent $\hat{c}=c_0-c_1-c_3$ and $\vec{N}(\alpha) = [0,k_{\alpha 2}\alpha^2+k_{\alpha 3}\alpha^3,0]^T$. The meaning of the parameters and their values used are given in Table~\ref{t1}. The bifurcation parameter $\mu$ represents here the wind velocity. $h$ and $\alpha$ stand for the heave displacement and the pitch angle, respectively.

\begin{table}
  \centering
  \begin{tabular}{ccl}
    \toprule 
    Parameter & Value &Description \\
    \midrule
    $U$            & $0$--$25$ & Airspeed (m/s)\\
    $b$            & $0.15$ & Wing semi-chord (m) \\
    $a$            & $-0.5$ & Position of elastic axis relative to the semi-chord (nd) \\
    $\rho$         & $1.204$ & Air density (kg/m\textsuperscript3) \\
    $m_w$          & $5.3$ & Mass of the wing (kg) \\
    $m_T$          & $16.9$ & Mass of wing and support (kg) \\
    $I_\alpha$     & $0.1726$  & Wing moment of inertia about elastic axis (kg\,m\textsuperscript2) \\
    $c_\alpha$     & $0.5628$ & Pitch linear damping coefficient (kg\,m\textsuperscript2/s) \\
    $c_h$          & $15.443$& Heave linear damping coefficient (kg/s) \\
    $k_\alpha$     & $54.1162$& Pitch linear stiffness (N/rad) \\
    $k_{\alpha 2}$ & $751.6$ & Pitch quadratic nonlinear stiffness (N/rad\textsuperscript2) \\
    $k_{\alpha 3}$ & $5006.7$& Pitch cubic nonlinear stiffness (N/rad\textsuperscript3) \\
    $k_{h}$        & $3529.4$& Heave linear stiffness (N/m) \\
    $x_\alpha$     & $0.234$ & Distance between center of gravity and elastic axis (nd) \\
    $c_{0,...,4}$  & $(1,\, 0.1650,\, 0.0455,\, 0.335,\, 0.3)$ & Aeroelastic coefficients\\
    \bottomrule
  \end{tabular}
  \caption{Descriptions of the parameters of \cref{eq:flutter2} and their values where applicable. Non-dimensional units are indicated by `nd'.}
  \label{t1}
\end{table}

This system has a subcritical Hopf bifurcation followed by saddle-node bifurcation of periodic orbits. The training data set includes four LCOs recorded on the stable branch, and four LCOs recorded on the unstable branch. As further discussed in \cref{sec:ExpVal}, in an experiment, stable and unstable LCOs can be directly measured using control-based continuation~\cite{RSPA-paper,lee2020reduced}. Here, the unstable LCOs were obtained by simulating the model~\cref{eq:flutter2} under proportional-derivative feedback control to reproduce the process followed for the experimental tests in \cref{sec:ExpVal}. Each time series is recorded for one second with a sampling time of 0.001$s$.

\begin{figure}[hbt!]
  \centering
    \includegraphics[width=0.65\linewidth]{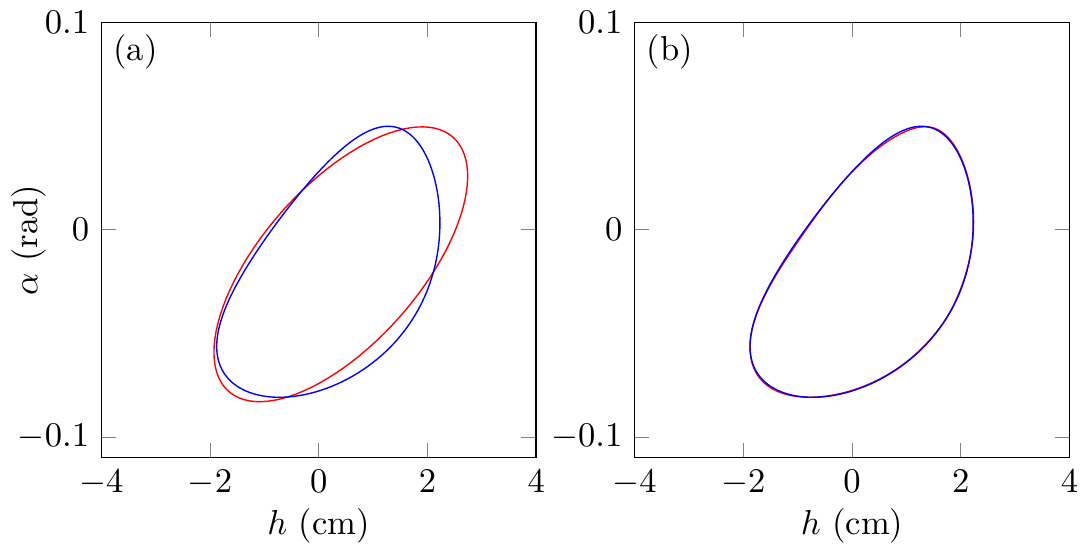}
\caption{Comparison between the phase portraits of the aeroelastic model (\textcolor{blue}{$\boldsymbol{-}$}) and the hybrid M/ML model (\textcolor{red}{$\boldsymbol{-}$}) for a stable LCO at $\mu=15.5 \; m/s$. Coordinate transformation (a) without and (b) with the neural network.}\label{LT}
\end{figure}

For this example, $\textrm{NN}_{\Theta_{\vec U}}$ has two inputs, two hidden layers each with 21 neurons and hyperbolic tangent activation functions, and three linear outputs. 400 iterations in ADAM with a learning rate of $0.01$ were necessary to estimate the NN parameters $\Theta_{\vec U}$. \cref{LT}(left) shows that the linear transformation allows the coordinate transformation to capture the overall orientation and size of the LCO, and \cref{LT}(right) shows that $\textrm{NN}_{\Theta_{\vec U}}$ further improves the accuracy of this coordinate transformation, leading to an excellent visual agreement between the LCO of the hybrid M/ML and reference models.

\begin{figure}[hbt!]
  \centering
  \subfloat{\includegraphics[width=0.7\textwidth]{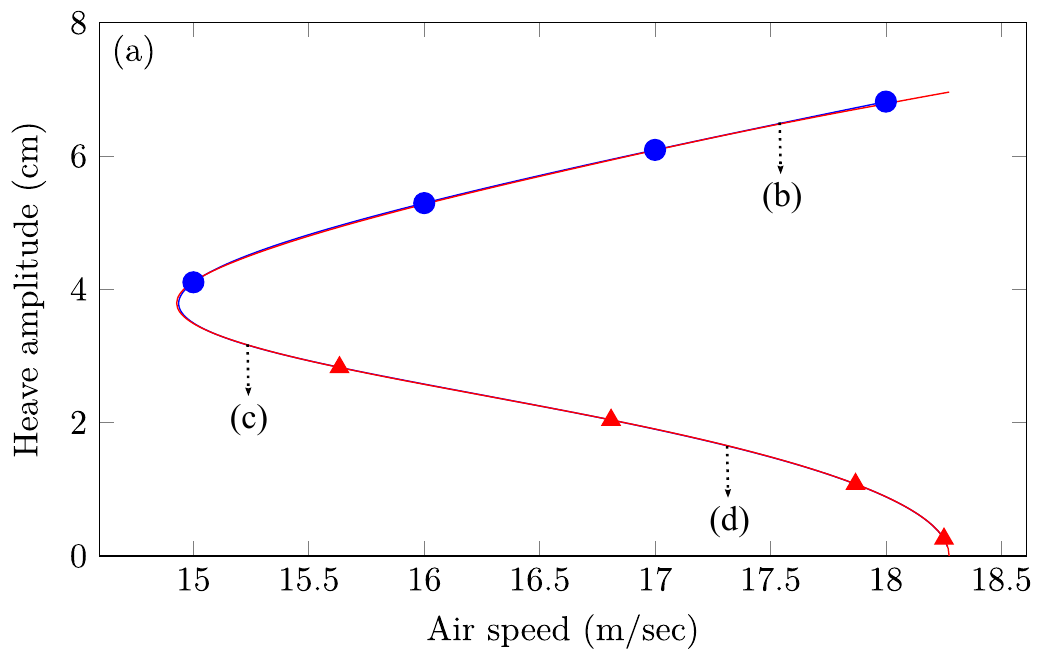}}
  \vspace*{-0.4em}
  \subfloat{\includegraphics[width=1.0\textwidth]{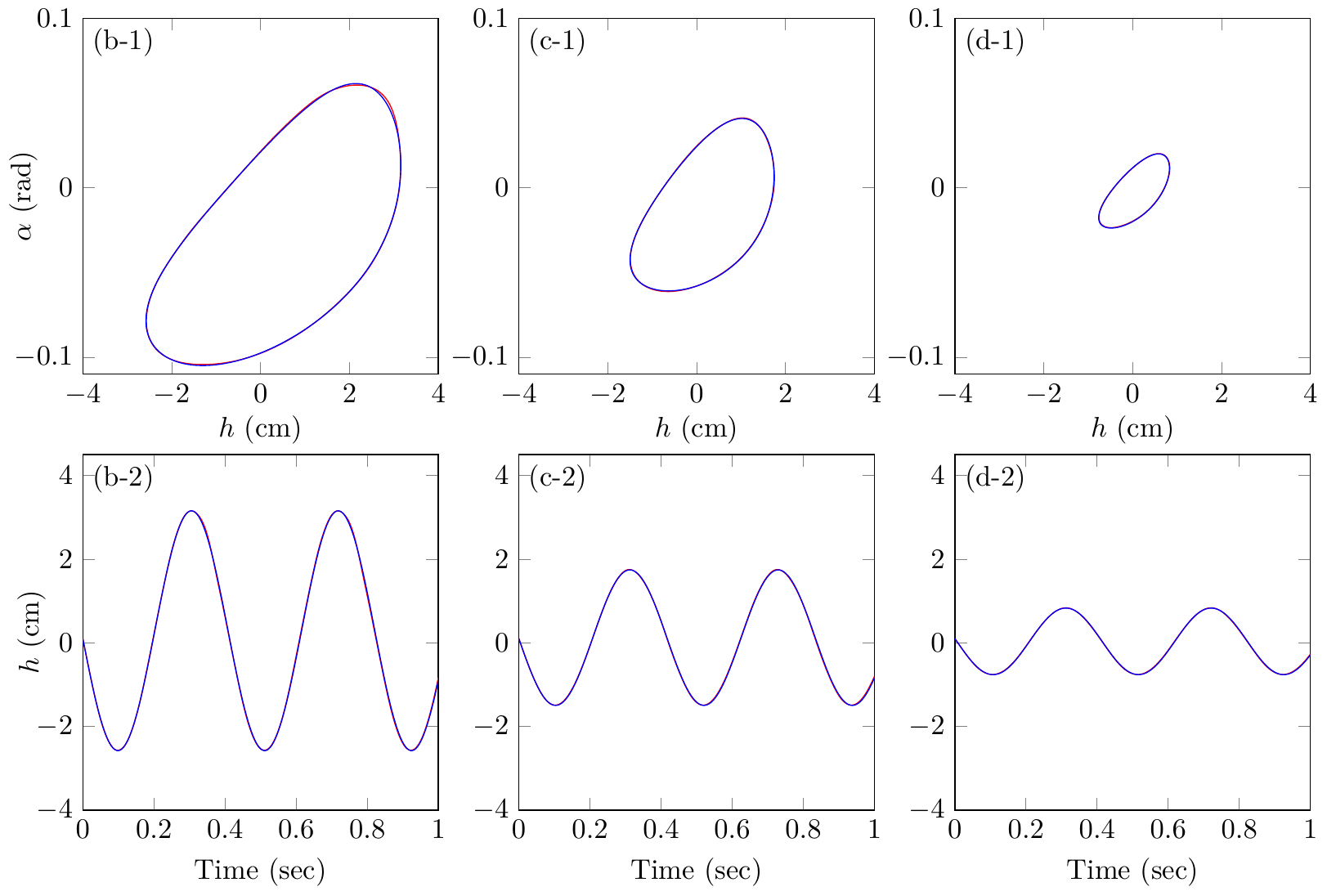}}
  \caption{Comparison between the aeroelastic model (\textcolor{blue}{$\boldsymbol{-}$}) and the hybrid M/ML model (\textcolor{red}{$\boldsymbol{-}$}). (a) Bifurcation diagram (\textcolor{blue}{$\boldsymbol{\bullet}$}) are stable LCOs and (\textcolor{red}{$\boldsymbol{\blacktriangle}$}) are unstable LCOs used for model training) (b-1--d-1) Phase portraits and (b-2--d-2) time-series at the locations reported on the bifurcation diagram.}
  \label{fig:flutter_all}
\end{figure}

\cref{fig:flutter_all}(a) shows there is an excellent agreement between the bifurcation diagrams computed from the reference and hybrid M/ML models. The bifurcation diagram of the trained model was computed by transforming 100 equi-spaced points on the periodic solutions of the normal-form model using $\vec U$. The identified values of the Hopf bifurcation point, $\mu_0$, and the saddle-node bifurcation point, $a_2$, are $18.28 \; m/s $ and $3.64$, respectively. This is in excellent agreement with the model values $\mu_0=18.28 \; m/s$ and $a_2=3.65$. The phase portrait of the trained model shows good agreement with the model for both stable and unstable LCOs (Figs.~\ref{fig:flutter_all}(b-1)--(d-1)).

The oscillation speed $\Omega$ is modelled using \cref{eq:Om_V}. For this example, only the constant term in the Fourier expansion is kept such that $\Omega(\vec u,\mu)=\omega_0+\textrm{NN}_{\Theta_{\Omega}}(\vec u,\mu)$. The neural network consist of three inputs, two hidden layers with 31 neurons each and a hyperbolic tangent activation functions, and one linear output. 300 iterations in ADAM with a learning rate of 0.01 followed by 1000 BFGS iterations with a $0.001$ learning rate were necessary to train $\Omega$ by minimizing $\Xi_{\Omega}$. Figs.~\ref{fig:flutter_all}(b-2)--(d-2) show that the model captures the overall time series and frequency of the LCO for the range of wind velocities considered. 

\section{Experimental demonstration on an aeroelastic structure}\label{sec:ExpVal}
The method developed in this paper is now demonstrated on a physical aeroelastic system.

\subsection{Experimental set-up and data collection}\label{sec:ExpSetup}
\begin{figure}[hbt!]
\centering
\subfloat[]{\includegraphics[width=0.52\textwidth]{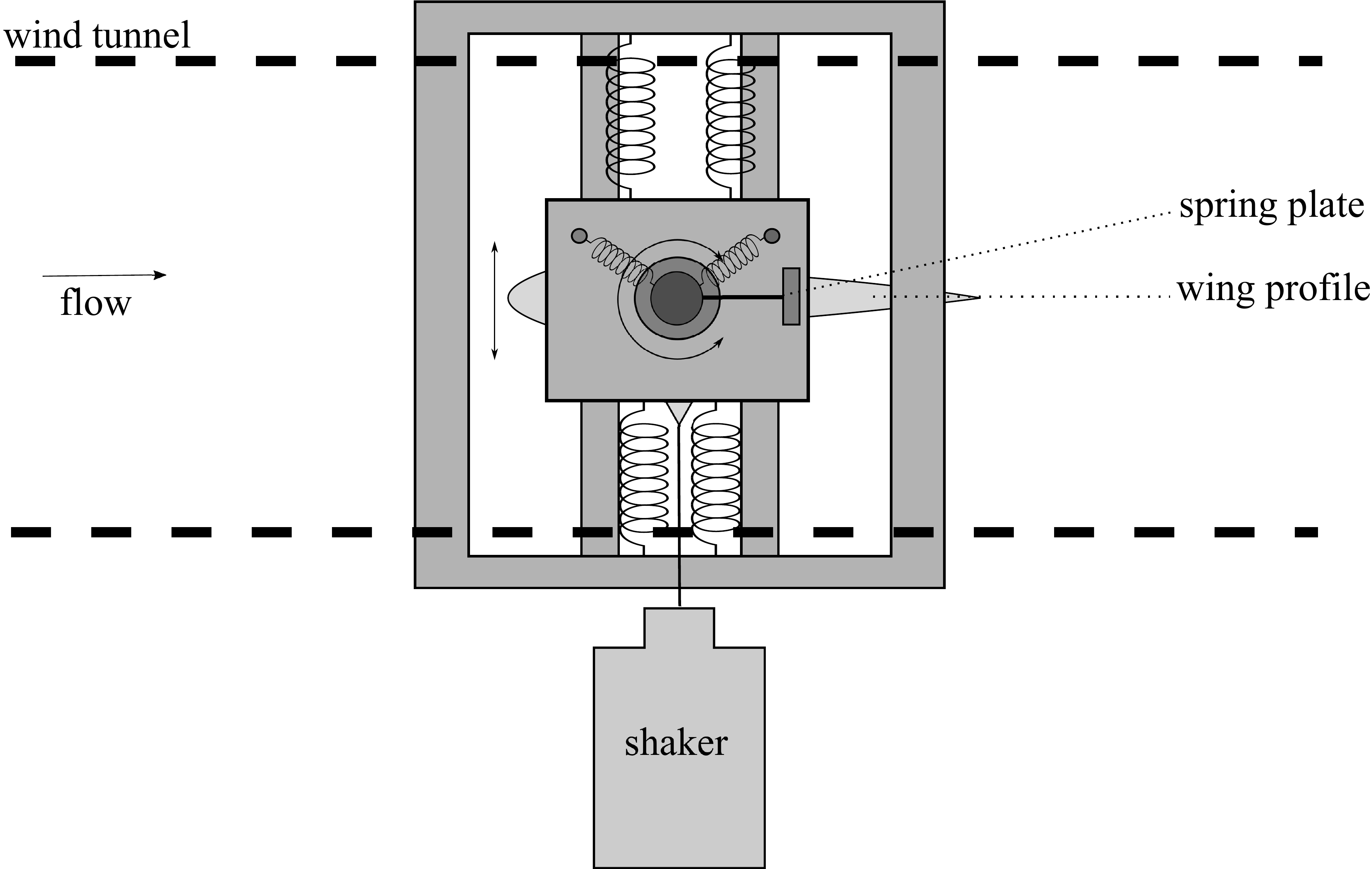}}
\hspace*{-0.3em}
\subfloat[]{\includegraphics[width=0.35\textwidth]{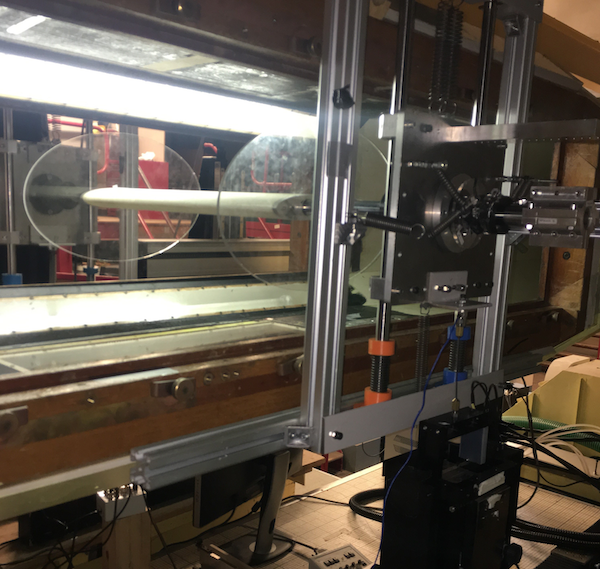}}
\caption{Aeroelastic rig. (a) Schematic. (b) Picture of the rig fitted to the University of Bristol's low-turbulence wind tunnel.}
\label{fig:rig}
\end{figure}

The rig is shown in Figure~\ref{fig:rig}. It comprises a NACA-0015 wing profile rigidly attached to a stainless steel shaft, supported at both ends by rotational bearings mounted on supporting plates that are constrained to move vertically by a linear bearing system. The structure has two mechanical degrees of freedom: one in pitch (rotational motion) and one in heave (vertical motion). In the heave direction, linear springs are connected between the supporting plates and the outer frame. In the pitch direction, torsional springs are connected between the shaft and the supporting plates. Both sets of springs provide approximately linear restoring forces in their respective directions. In the pitch direction there are additional leaf springs connecting the shaft and the supporting plates; these leaf springs provide a hardening nonlinearity, mimicking potential interface effects at the root of the aerofoil. The dimensions of the flutter rig are such that the wing profile fits in the principal section of the University of Bristol's low-turbulence wind tunnel; the supporting plates and outer frame lie outside it (see Figure~\cref{fig:rig}(b)). The reader is referred to \cite{RSPA-paper} for further details about the dimensions of the system.

Control-based continuation was exploited to measure the stable and unstable LCOs of the system directly during the wind tunnel tests. The control forces are applied in the heave direction by use of an APS 113 electro-seis long-stroke electrodynamic shaker connected by a flexible stinger to one of the supporting plates. The experiment is instrumented with an Omron ZX1-LD300 laser displacement sensor to measure the heave motion, and an RLS AksIM 18\,bit absolute magnetic encoder fitted on the shaft to capture the pitch motion. The wind speed was directly provided by the wind tunnel control system. Real-time control and data acquisition is performed using a Beaglebone Black single-board computer equipped with an analogue IO cape (18\,bit ADC and 16\,bit DAC) operating at a sample rate of 5\,kHz~\cite{CBC_hardware}. The reader is referred to~\cite{RSPA-paper} for further details about the experimental set-up and the CBC method.

\subsection{The hybrid M/ML model}
The dynamics of this aeroelastic system is characterised by a subcritical Hopf bifurcation followed by a saddle-node bifurcation of cycles. The mechanistic model used within the hybrid model is therefore the one presented in \cref{NF-model} and already used in \cref{sec:NumVal_aero}. The training data sets includes four stable LCOs measured at $\mu= (14.9,\;15.6,\;16.5,\;17.3) \; m/s$, and three unstable LCOs measured at $\mu = (14.9,\;15.6,\;16.5) \; m/s$. Time series includes 6000 time points per LCOs, which represents approximately 12 oscillation periods. To train the oscillation speed, time series were down sampled to 1000 samples to reduce the computational cost of the training.

For the coordinate map $\vec U_{12}$, a neural network $\textrm{NN}_{\Theta_{\vec U}}$ with three inputs, two outputs and two hidden layers each with 11 neurons and hyperbolic tangent activation functions was used. A first 1000 iterations with ADAM using a $0.01$ learning rate followed by 3000 BFGS iterations with a $0.0001$ learning rate were necessary to minimize $\Xi_{\vec U}$ and find $\mu_0,\; a_2$ and the network parameters $\Theta_{\vec U}$.

\begin{figure}[t]
    \centering
      \includegraphics[width=0.7\linewidth]{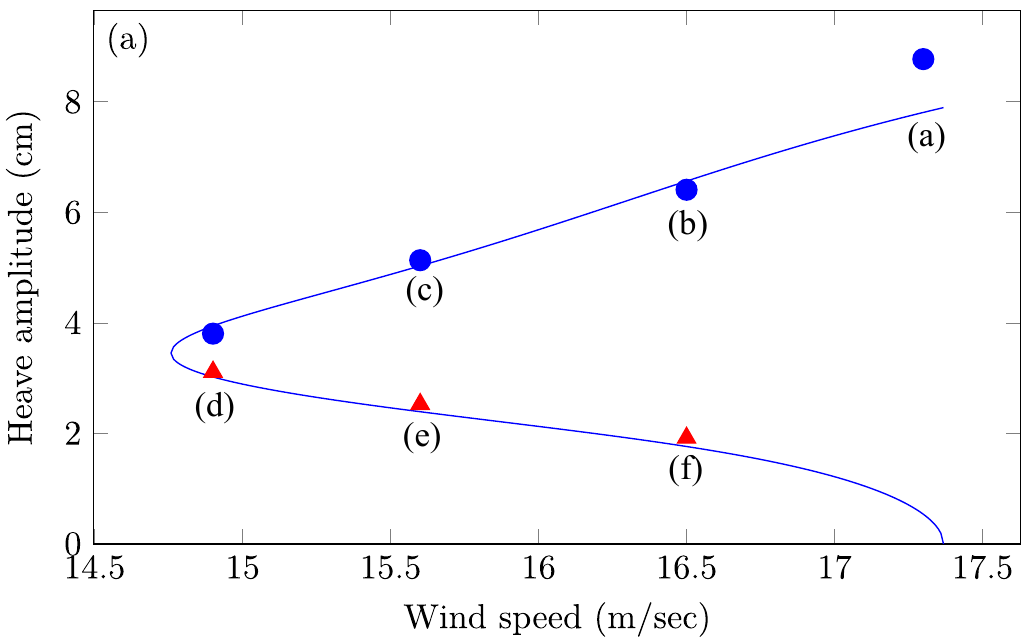}
  \caption{Comparison between the bifurcation diagram obtained from the hybrid M/ML model with increased model accuracy  (\textcolor{blue}{$\boldsymbol{-}$}) and the stable (\textcolor{blue}{$\boldsymbol{\bullet}$}) and unstable (\textcolor{red}{$\boldsymbol{\blacktriangle}$}) LCOs used for model training. Labels (a)-(f) denote the corresponding phase portraits and time series plots in \cref{pp_exp} and~\cref{ts_exp}, respectively.}
    \label{bd_exp}
  \end{figure}

\begin{figure}[h!]
    \centering
      \includegraphics[width=0.9\linewidth]{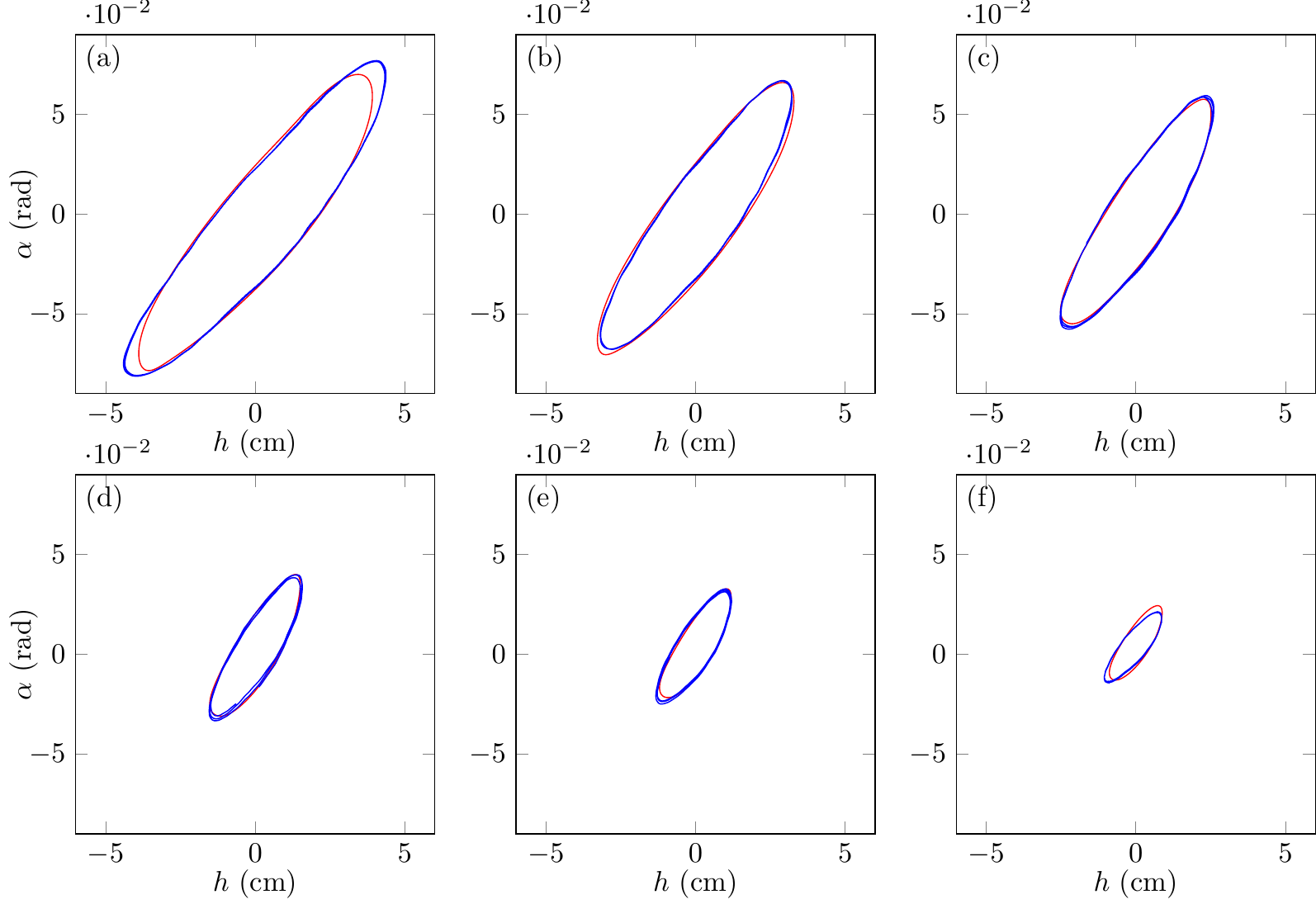}
  \caption{Comparison of phase portraits between the measured LCOs (\textcolor{blue}{$\boldsymbol{-}$}) and the hybrid M/ML model (\textcolor{red}{$\boldsymbol{-}$}). (a) stable LCO at wind speed 17.3 m/sec, (b) stable LCO at wind speed 16.5 m/sec, (c) stable LCO at wind speed 15.6 m/sec, (d) unstable LCO at wind speed 14.9 m/sec, (e) unstable LCO at wind speed 15.6 m/sec and (f) unstable LCO at wind speed 16.5 m/sec.}
    \label{pp_exp}
  \end{figure}

\cref{bd_exp} compares the bifurcation diagram of the hybrid M/ML model with the LCOs measured experimentally. A qualitatively good agreement with the data is obtained despite the limited number of LCOs used for model training. The Hopf bifurcation point is estimated at $\mu_0=17.67 \; m/s$ and the saddle-node bifurcation point at $14.66 \; m/s$. Overall, the trained model accurately predicts the phase portrait of the stable and unstable LCOs, as shown in~\cref{pp_exp}. In the phase portraits, the line associated with the experimentally measured LCOs appears thicker than the one from model predictions. This is an illusion that comes from the presence of multiple oscillation periods in the recorded data and the unavoidable differences that exist between periods due to the presence of noise in the measurements.

\begin{figure}[h!]
    \centering
      \includegraphics[width=0.9\linewidth]{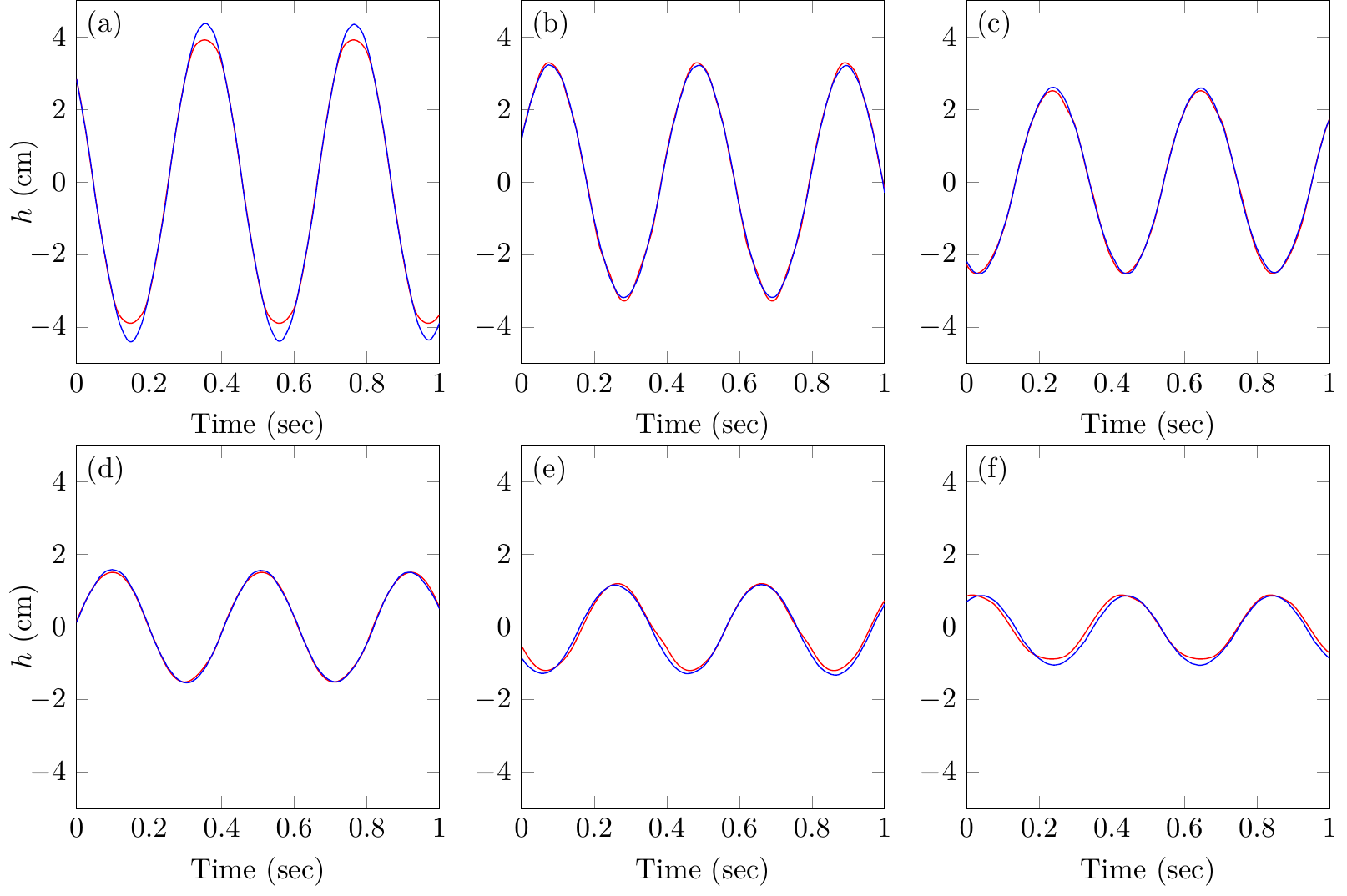}
  \caption{Comparison of heave time series between the measured LCOs (\textcolor{blue}{$\boldsymbol{-}$}) and the hybrid M/ML model (\textcolor{red}{$\boldsymbol{-}$}). (a) stable LCO at wind speed 17.3 m/sec, (b) stable LCO at wind speed 16.5 m/sec, (c) stable LCO at wind speed 15.6 m/sec, (d) unstable LCO at wind speed 14.9 m/sec, (e) unstable LCO at wind speed 15.6 m/sec and (f) unstable LCO at wind speed 16.5 m/sec.}
    \label{ts_exp}
  \end{figure}

For the identification of $\Omega$, a similar model to the one used in \cref{sec:NumVal_aero} is considered. The neural network $\textrm{NN}_{\Theta_{\Omega}}(\vec u,\mu)$ includes three inputs, a single linear output and two hidden layers with 21 neurons each and hyperbolic tangent activation functions. 500 iterations in ADAM with a learning rate of 0.01 followed by 400 BFGS iterations with a $0.001$ learning rate were necessary to minimize $\Xi_{\Omega}$. The time series presented in \cref{ts_exp} show that the model captures the frequency of the measured LCOs. The amplitude error visible in the bifurcation diagram is also clearly visible in the time series.

\begin{figure}[h!]
	\centering
    \includegraphics[width=0.7\linewidth]{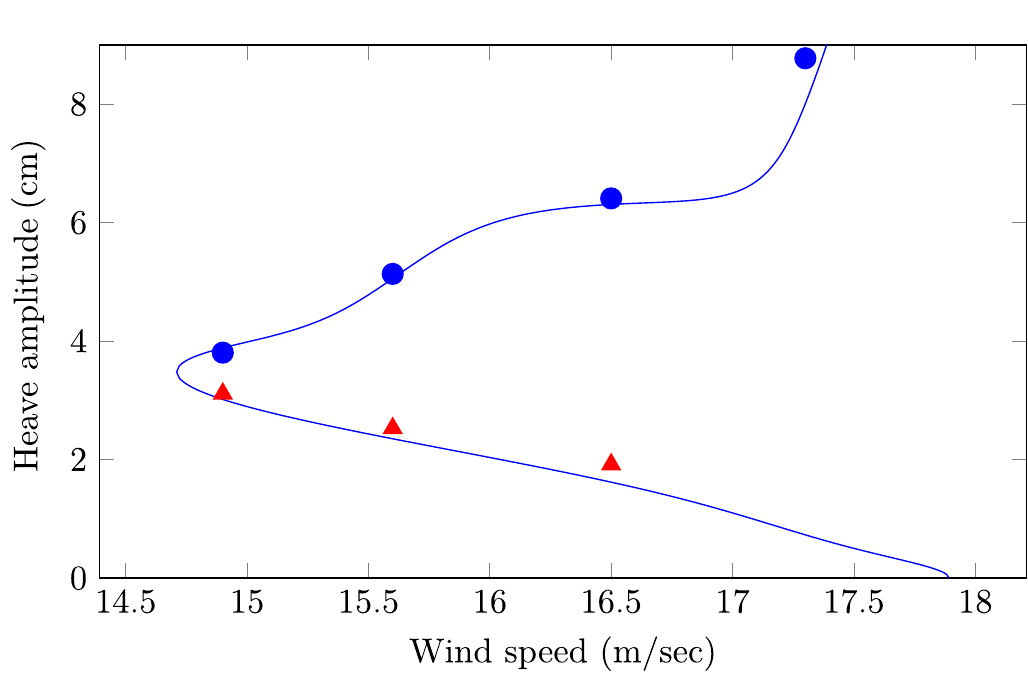}
\caption{Illustration of model overfitting. Bifurcation diagram obtained from the hybrid M/ML model (\textcolor{blue}{$\boldsymbol{-}$}), stable (\textcolor{blue}{$\boldsymbol{\bullet}$}) and unstable (\textcolor{red}{$\boldsymbol{\blacktriangle}$}) LCOs used for model training.}
  \label{bd_of}
\end{figure}

One potential issue with ML model training is overfitting. This is illustrated in~\cref{bd_of} where the model was trained using different initial parameter values and different hyperparameters (number of iterations and learning rate). While the overall model prediction error is small at the data points, the model presents a large variability between them which is symptomatic of overfitting. To reduce overfitting and obtain the results presented in Figures~\ref{bd_exp}~--~\ref{ts_exp}, some hyper-parameters, such as the number of iterations, were manually tuned. Approaches that promote parameter sparsity~\cite{louizos2017learning} or a formal optimisation of the hyperparameters~\cite{burden2008bayesian} were not carried out due to the associated computational costs and the overall lack of data. Note that, the use of physics --- here, through the model structure~\eqref{NF-model} and the use of linear transformations in $\vec U_{12}$ --- can also be viewed as a regularisation techniques that reduces ML model complexity and hence helps in reducing overfitting. 

\begin{figure}[h!]
    \centering
    \subfloat{\includegraphics[width=0.7\textwidth]{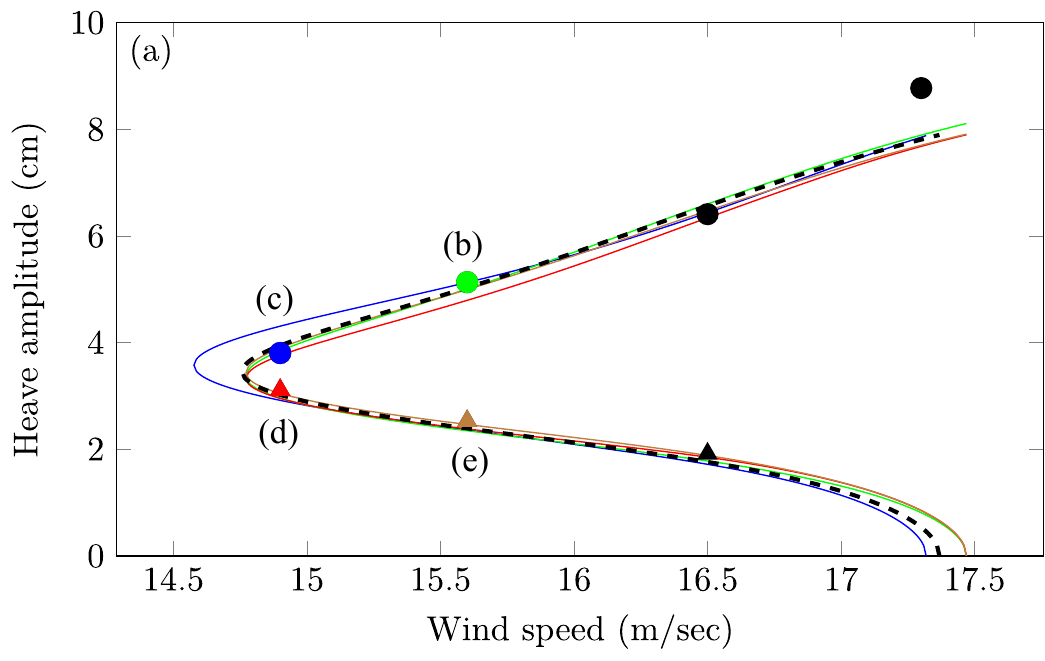}}
    \vspace*{-0.4em}
    \subfloat{\includegraphics[width=1.0\textwidth]{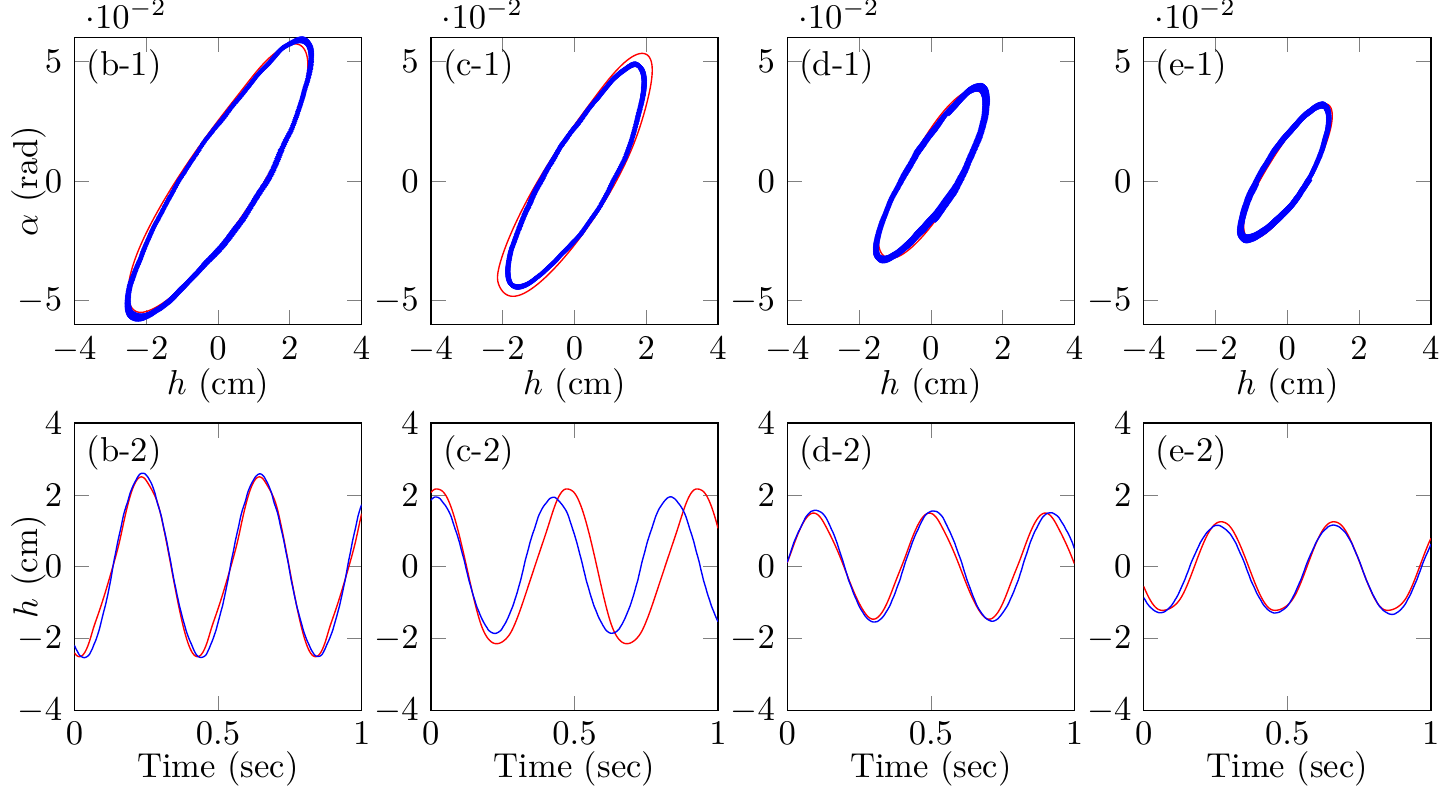}}
    \caption{Effect of excluding one data point from the training data set on the hybrid M/ML model accuracy. (a) Bifurcation diagrams obtained when removing the data point of the same colour. For instance, the blue bifurcation curve is obtained with the model trained with a data set excluding the blue point. (\textcolor{black}{$\boldsymbol{--}$}) Bifurcation diagram obtained from a model trained with all measured data. (b-1) -- (e-1) Prediction of phase portrait at the excluded data point. (b-2) -- (e-2) Prediction of time series at the excluded data point.} 
    \label{fig:exp_val}
  \end{figure}

To assess the robustness of the identified model with respect to the training data, the model training was also performed with four different data sets, each with one of the LCO data points removed. This approach is inspired by the leave-one-out cross-validation technique and chosen due to the small number of data points available in the parameter space. \cref{fig:exp_val} shows the bifurcation diagrams obtained after removing the different data points. The colour of the bifurcation curve matches the colour of the data point that was removed from the training data set. The dashed black bifurcation curve was obtained by including all the data points in the training set. The phase portraits and time series shown in~\cref{fig:exp_val} illustrate the performance of the model at the removed data point. While most bifurcation curves appear similar, removing the stable LCO in blue appears to have a significant influence on the location of the saddle-node bifurcation and more generally on the bifurcation curve in that area. This also affects the quality of the oscillation speed model $\Omega$, which is unable to capture the LCO oscillation frequency adequately (see \cref{fig:exp_val}(c-2)). Overall this results suggest that sufficient training data near bifurcation points (saddle-node and Hopf) is needed to build a robust model. 

\begin{figure}[hbt!]
	\centering
    \includegraphics[width=0.9\linewidth]{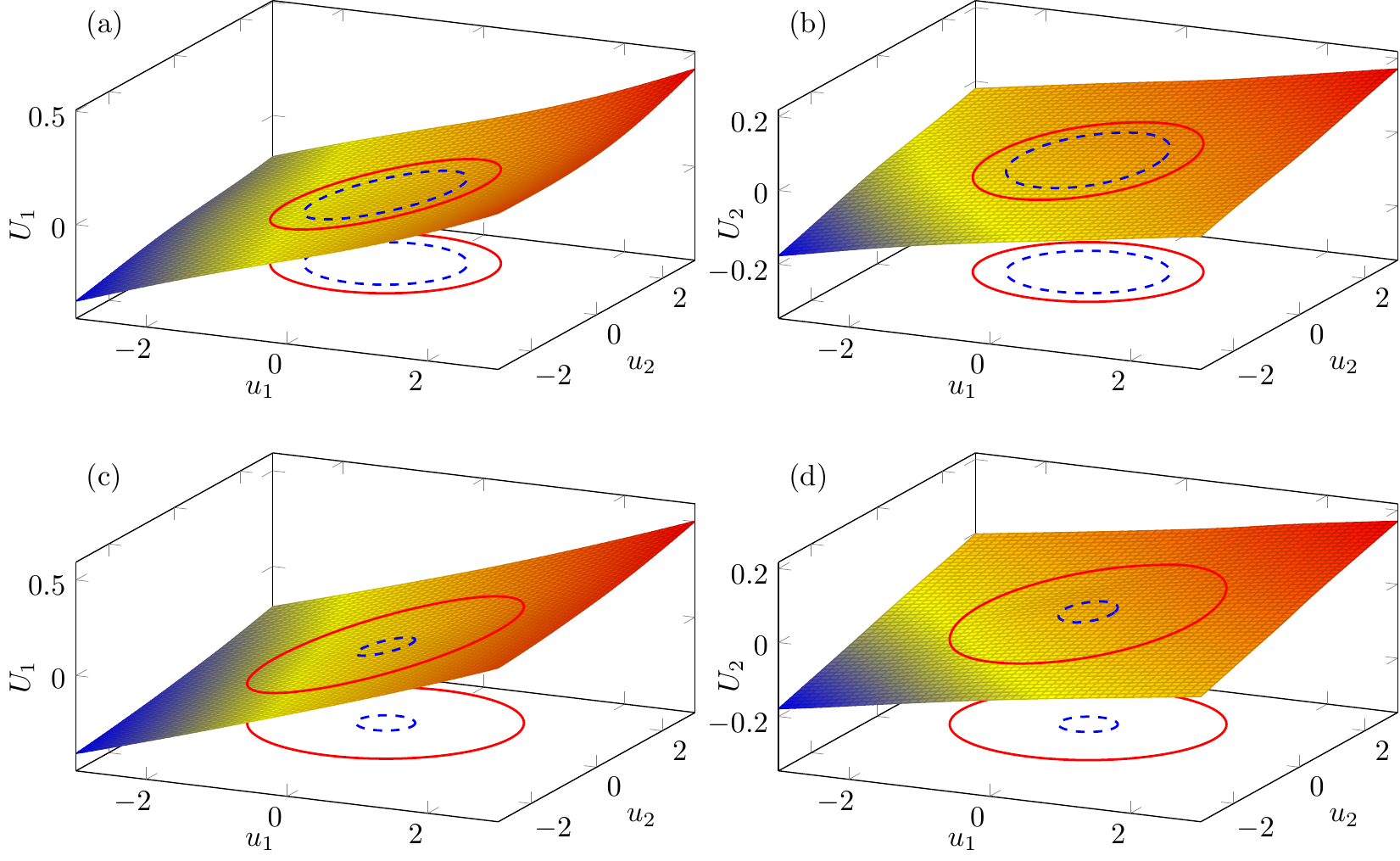}
\caption{Mapping $\vec U_{12}(\vec u,\mu)=(U_1(\vec u,\mu),U_2(\vec u,\mu))$ between normal form and physical coordinates. (a) $U_1$ at $\mu$=15.0 $m/sec$, (b) $U_2$ at $\mu$=15.0 $m/sec$, (c) $U_1$ at $\mu$=17.0 $m/sec$ and (d) $U_2$ at $\mu$=17.0 $m/sec$.}
\label{trans}
\end{figure}

\cref{trans} presents the trained mapping $\vec U_{12}$. \cref{trans}(a, b) show the transformation $\vec U_{12}$ for $\mu=$ 15.0, and \cref{trans}(c, d) show the transformation $\vec U_{12}$ for $\mu=$ 17.5 m/s. The blue dotted lines and the red solid lines correspond to the unstable and stable LCOs, respectively. The coordinate transformations are smooth transformation and locally invertible. The visible curvature shows that the transformations are also nonlinear. As discussed in Section~\ref{sec:map_func}, the presence of an initial linear coordinate transformation in~\eqref{g} was essential. A NN alone was unable to produce topologically equivalent closed curves and obtain a locally invertible transformation near the bifurcation point. Models with significant overfitting, such as the one in~\cref{bd_of}, were also found to result in poorly- or even non-invertible transformations.

\FloatBarrier

\section{Conclusions}
In this paper, we have proposed a new hybrid modelling approach for physical systems with a Hopf bifurcation. At its core, it uses a mechanistic model, in the form of a normal-form-like model, to capture the phenomenology of the physical system. A data-driven, machine-learnt model is then used to map the mechanistic model predictions onto the measured data. Our method was first demonstrated with numerical data collected on a Van der Pol oscillator and an aeroelastic model, and then with experimental data collected on aeroelastic rig during wind tunnel tests. The hybrid mechanistic/machine-learnt models obtained with our method were shown to quantitatively capture the bifurcation diagrams of the different systems as well as their time evolution, even in the presence of multiple time-scales and noise. 

The proposed method has several advantages such as being data-driven whilst also being able to work with a limited number of measured states and data. It also requires only the knowledge of the bifurcation structure of the system and is thus applicable to  any system with a Hopf bifurcation. Replacing the current mechanistic model by a more detailed model of the physical system could improve the accuracy of the obtained model while also simplifying the training of the mapping. However, this would come at the cost of having a more application-specific modelling methodology. A more systematic approach to handle overfitting during ML model training is also required to make the method more systematically applicable. As the model is derived using LCO data only, the obtained hybrid mechanistic/machine-learnt model reproduces the long-term behaviour of the physical system but is generally unable to accurately capture its transient dynamics. Future work should look at including transient data in the modelling approach. Finally, the principles of the proposed method are extremely general and could therefore be applied to systems with other types of bifurcations by changing the mechanistic model.

The hybrid mechanistic/machine-learnt models obtained with our method can have multiple uses. The machine-learnt part of the model could be exploited to improve understanding about the physical system of interest and provide new insights into the derivation of more accurate mechanistic models. The ability of the proposed models to reproduce the unstable part of the bifurcation diagram has also the potential help analyse stability boundaries and basins of attraction of physical systems, or help in reducing the cost of experimental methods such as CBC. Finally, the hybrid mechanistic/machine-learnt models developed here establish a rigorous framework to combine known physics with experimental and/or operational data and could therefore be used as digital twins for nonlinear systems. 

\section*{Statements}
\subsection*{Funding}
This work was supported by a PhD Scholarship from the University of Bristol, the EPSRC (EP/K032738/1) and the Royal Academy of Engineering (RF1516/15/11).

\subsection*{Data Availability}
The datasets generated during and/or analysed during the current study are available in \url{https://github.com/Kyounghyunlee/ML_Hopf}.

\bibliography{ref}

\begin{thebibliography}{10}

\bibitem{dimitriadis2017introduction}
Dimitriadis G.
\newblock Introduction to Nonlinear Aeroelasticity.
\newblock John Wiley \& Sons; 2017.

\bibitem{beregi2019bifurcation}
Beregi S, Takacs D, Stepan G.
\newblock Bifurcation analysis of wheel shimmy with non-smooth effects and time
  delay in the tyre--ground contact.
\newblock Nonlinear Dynamics. 2019;98(1):841-58.

\bibitem{kalmar2001subcritical}
Kalm{\'a}r-Nagy T, St{\'e}p{\'a}n G, Moon FC.
\newblock Subcritical Hopf bifurcation in the delay equation model for machine
  tool vibrations.
\newblock Nonlinear Dynamics. 2001;26(2):121-42.

\bibitem{adimy2005stability}
Adimy M, Crauste F, Ruan S.
\newblock Stability and Hopf bifurcation in a mathematical model of pluripotent
  stem cell dynamics.
\newblock Nonlinear Analysis: Real World Applications. 2005;6(4):651-70.

\bibitem{guo2003hopf}
Guo S, Huang L.
\newblock Hopf bifurcating periodic orbits in a ring of neurons with delays.
\newblock Physica D: Nonlinear Phenomena. 2003;183(1-2):19-44.

\bibitem{kuznetsov2013elements}
Kuznetsov YA.
\newblock Elements of applied bifurcation theory. vol. 112.
\newblock Springer Science \& Business Media; 2013.

\bibitem{winkler2017performance}
Winkler DA, Le TC.
\newblock Performance of deep and shallow neural networks, the universal
  approximation theorem, activity cliffs, and QSAR.
\newblock Molecular informatics. 2017;36(1-2):1600118.

\bibitem{wang2021understanding}
Wang S, Teng Y, Perdikaris P.
\newblock Understanding and mitigating gradient flow pathologies in
  physics-informed neural networks.
\newblock SIAM Journal on Scientific Computing. 2021;43(5):A3055-81.

\bibitem{rasmussen2003gaussian}
Rasmussen CE.
\newblock Gaussian processes in machine learning.
\newblock In: Summer school on machine learning. Springer; 2003. p. 63-71.

\bibitem{kim2019dpm}
Kim J, Lee K, Lee D, Jin SY, Park N.
\newblock DPM: A Novel Training Method for Physics-Informed Neural Networks in
  Extrapolation.
\newblock Comput Phys. 2019;378:686-707.

\bibitem{raissi2019physics}
Raissi M, Perdikaris P, Karniadakis GE.
\newblock Physics-informed neural networks: A deep learning framework for
  solving forward and inverse problems involving nonlinear partial differential
  equations.
\newblock Journal of Computational Physics. 2019;378:686-707.

\bibitem{raissi2017machine}
Raissi M, Perdikaris P, Karniadakis GE.
\newblock Machine learning of linear differential equations using Gaussian
  processes.
\newblock Journal of Computational Physics. 2017;348:683-93.

\bibitem{raissi2017inferring}
Raissi M, Perdikaris P, Karniadakis GE.
\newblock Inferring solutions of differential equations using noisy
  multi-fidelity data.
\newblock Journal of Computational Physics. 2017;335:736-46.

\bibitem{rackauckas2020universal}
Rackauckas C, Ma Y, Martensen J, Warner C, Zubov K, Supekar R, et~al.
\newblock Universal differential equations for scientific machine learning.
\newblock arXiv preprint arXiv:200104385. 2020.

\bibitem{beregi2021using}
Beregi S, Barton DA, Rezgui D, Neild SA.
\newblock Using scientific machine learning for experimental bifurcation
  analysis of dynamic systems.
\newblock arXiv preprint arXiv:211011854. 2021.

\bibitem{sieber2008control}
Sieber J, Krauskopf B.
\newblock Control based bifurcation analysis for experiments.
\newblock Nonlinear Dynamics. 2008;51(3):365-77.

\bibitem{renson2019application}
Renson L, Shaw A, Barton D, Neild S.
\newblock Application of control-based continuation to a nonlinear structure
  with harmonically coupled modes.
\newblock Mechanical Systems and Signal Processing. 2019;120:449-64.

\bibitem{renson2016robust}
Renson L, Gonzalez-Buelga A, Barton D, Neild S.
\newblock Robust identification of backbone curves using control-based
  continuation.
\newblock Journal of Sound and Vibration. 2016;367:145-58.

\bibitem{barton2017control}
Barton DA.
\newblock Control-based continuation: Bifurcation and stability analysis for
  physical experiments.
\newblock Mechanical Systems and Signal Processing. 2017;84:54-64.

\bibitem{barton2013systematic}
Barton DA, Sieber J.
\newblock Systematic experimental exploration of bifurcations with noninvasive
  control.
\newblock Physical Review E. 2013;87(5):052916.

\bibitem{deCesare22}
de~Cesare I, Salzano D, di~Bernardo M, Renson L, Marucci L.
\newblock Control-Based Continuation: A New Approach to Prototype Synthetic
  Gene Networks.
\newblock ACS Synthetic Biology, in press. 0;0(0):null.

\bibitem{beregi2020improving}
Beregi S, Barton DA, Rezgui D, Neild SA.
\newblock Improving robustness to noise of nonlinear parameter identification
  using control-based continuation.
\newblock arXiv preprint arXiv:200111008. 2020.

\bibitem{carr2012applications}
Carr J.
\newblock Applications of centre manifold theory. vol.~35.
\newblock Springer Science \& Business Media; 2012.

\bibitem{rosenfeld1976digital}
Rosenfeld A.
\newblock Digital picture processing.
\newblock Academic press; 1976.

\bibitem{zhang2004review}
Zhang D, Lu G.
\newblock Review of shape representation and description techniques.
\newblock Pattern recognition. 2004;37(1):1-19.

\bibitem{zahn1972fourier}
Zahn CT, Roskies RZ.
\newblock Fourier descriptors for plane closed curves.
\newblock IEEE Transactions on computers. 1972;100(3):269-81.

\bibitem{penrose1955generalized}
Penrose R.
\newblock A generalized inverse for matrices.
\newblock In: Mathematical proceedings of the Cambridge philosophical society.
  vol.~51. Cambridge University Press; 1955. p. 406-13.

\bibitem{kanagawa2018gaussian}
Kanagawa M, Hennig P, Sejdinovic D, Sriperumbudur BK.
\newblock Gaussian processes and kernel methods: A review on connections and
  equivalences.
\newblock arXiv preprint arXiv:180702582. 2018.

\bibitem{lin2018resnet}
Lin H, Jegelka S.
\newblock Resnet with one-neuron hidden layers is a universal approximator.
\newblock In: Advances in neural information processing systems; 2018. p.
  6169-78.

\bibitem{paszke2019pytorch}
Paszke A, Gross S, Massa F, Lerer A, Bradbury J, Chanan G, et~al.
\newblock Pytorch: An imperative style, high-performance deep learning library.
\newblock In: Advances in neural information processing systems; 2019. p.
  8026-37.

\bibitem{innes2018flux}
Innes M.
\newblock Flux: Elegant machine learning with Julia.
\newblock Journal of Open Source Software. 2018;3(25):602.

\bibitem{da2014method}
Da K.
\newblock A method for stochastic optimization.
\newblock arXiv preprint arXiv:14126980. 2014.

\bibitem{liu1989limited}
Liu DC, Nocedal J.
\newblock On the limited memory BFGS method for large scale optimization.
\newblock Mathematical programming. 1989;45(1-3):503-28.

\bibitem{zhang2017discrete}
Zhang H, Abhyankar S, Constantinescu E, Anitescu M.
\newblock Discrete adjoint sensitivity analysis of hybrid dynamical systems
  with switching.
\newblock IEEE Transactions on Circuits and Systems I: Regular Papers.
  2017;64(5):1247-59.

\bibitem{lauss2018discrete}
Lau{\ss} T, Oberpeilsteiner S, Steiner W, Nachbagauer K.
\newblock The discrete adjoint method for parameter identification in multibody
  system dynamics.
\newblock Multibody system dynamics. 2018;42(4):397-410.

\bibitem{ruder2016overview}
Ruder S.
\newblock An overview of gradient descent optimization algorithms.
\newblock arXiv preprint arXiv:160904747. 2016.

\bibitem{fletcher2013practical}
Fletcher R.
\newblock Practical methods of optimization.
\newblock John Wiley \& Sons; 2013.

\bibitem{abdelkefi2013analytical}
Abdelkefi A, Vasconcellos R, Nayfeh AH, Hajj MR.
\newblock An analytical and experimental investigation into limit-cycle
  oscillations of an aeroelastic system.
\newblock Nonlinear Dynamics. 2013;71(1-2):159-73.

\bibitem{RSPA-paper}
Lee K, Tartaruga I, Rezgui D, Renson L, Neild SA, Barton DAW.
\newblock Analysis of self-excited flutter oscillations with control-based
  continuation.
\newblock Preprint. 2022.

\bibitem{lee2020reduced}
Lee K, Barton D, Renson L.
\newblock Reduced-order modelling of flutter oscillations using normal forms
  and scientific machine learning.
\newblock arXiv e-prints. 2020:arXiv-2011.

\bibitem{CBC_hardware}
Barton DAW. Real-time control hardware/software based on the BeagleBone Black;
  2015.
\newblock Available from: \url{http://github.com/~dawbarton/rtc}.

\bibitem{louizos2017learning}
Louizos C, Welling M, Kingma DP.
\newblock Learning sparse neural networks through $ L\_0 $ regularization.
\newblock arXiv preprint arXiv:171201312. 2017.

\bibitem{burden2008bayesian}
Burden F, Winkler D.
\newblock Bayesian regularization of neural networks.
\newblock Artificial neural networks. 2008:23-42.

\end{thebibliography}
\bibliographystyle{vancouver}

\end{document}